\RequirePackage{luatex85}
\documentclass[doublespace,12pt]{amsart}
\usepackage[numbers]{natbib}
\usepackage{bibunits}
\usepackage[english]{babel}
\usepackage[T1]{fontenc}
\usepackage[utf8]{inputenc}
\usepackage{amsmath,amssymb,amsthm}
\usepackage{textcomp}
\usepackage{stmaryrd}
\usepackage[sc]{mathpazo}
\usepackage[left=4cm,right=4cm,top=2cm,bottom=3cm]{geometry}
\usepackage{eulervm}
\usepackage[cal=dutchcal,scr=kp]{mathalpha}
\makeatletter
\def\mathalfa@frakscaled{s*[1.0]}
  \DeclareFontFamily{U}{esstixfrak}{\skewchar \font =45}
  \DeclareFontShape{U}{esstixfrak}{m}{n}{
    <-> \mathalfa@frakscaled esstixfrak}{}
  \DeclareMathAlphabet{\mathfrak}{U}{esstixfrak}{m}{n}
\makeatother

\usepackage{mathtools}
\usepackage{comment}

\usepackage{hyperref}
\usepackage{tikz}
\usetikzlibrary{knots,hobby,decorations.pathreplacing,shapes.geometric,calc,decorations.text,shapes.misc,calc,math,arrows,decorations.markings,patterns, matrix, intersections,cd}
\usepackage[nameinlink]{cleveref}

\definecolor{couleur3}{rgb}{0.0, 0.5, 1.0} 
\definecolor{couleur4}{rgb}{0.55, 0.71, 0.0} 
\definecolor{couleur1}{rgb}{0.87, 0.45, 1.0} 
\definecolor{couleur2}{rgb}{0.81, 0.06, 0.13}
\definecolor{couleur5}{rgb}{1, .64, 0.28}

\newcommand{\catname}[1]{\mathbf{#1}}

\newcommand{\sitename}[1]{\mathfrak{#1}}
\newcommand{\nom}[1]{\textsc{#1}}

\newcommand{\C}{\mathbb{C}}

\newcommand{\N}{\mathbb{N}}

\renewcommand{\P}{\mathbb{P}}
\newcommand{\Q}{\mathbb{Q}}
\newcommand{\R}{\mathbb{R}}
\renewcommand{\S}{\mathbb{S}}

\newcommand{\Z}{\mathbb{Z}}
\newcommand{\Acal}{\mathcal{A}}
\newcommand{\Bcal}{\mathcal{B}}

\newcommand{\Ical}{\mathcal{I}}

\newcommand{\Mcal}{\mathcal{M}}

\newcommand{\Pcal}{\mathcal{P}}
\newcommand{\pcal}{\mathcal{p}}

\newcommand{\Rcal}{\mathcal{R}}

\newcommand{\Ascr}{\mathscr{A}}
\newcommand{\Bscr}{\mathscr{B}}

\newcommand{\Fscr}{\mathscr{F}}

\newcommand{\Mscr}{\mathscr{M}}

\newcommand{\Sscr}{\mathscr{S}}
\newcommand{\Tscr}{\mathscr{T}}
\newcommand{\Uscr}{\mathscr{U}}

\newcommand{\Xscr}{\mathscr{X}}

\newcommand{\Sfrak}{\mathfrak{S}}

\renewcommand{\hbar}{\overline{h}}

\DeclareMathOperator{\Hom}{\mathrm{Hom}}

\DeclareMathOperator{\Vect}{\mathrm{Vect}}
\DeclareMathOperator{\Spec}{\mathrm{Spec}}

\DeclareMathOperator{\GL}{\mathrm{GL}}

\DeclareMathOperator{\im}{\mathrm{im}}
\DeclareMathOperator{\Cone}{\mathrm{Cone}}

\DeclareMathOperator{\colim}{\mathrm{colim}}

\DeclareMathOperator{\Poset}{\mathbf{Poset}}

\DeclareMathOperator{\Aut}{\mathrm{Aut}}
\DeclareMathOperator{\comb}{\mathrm{comb}}
\DeclareMathOperator{\Gr}{\mathrm{Gr}}

\DeclareMathOperator{\Card}{\mathrm{Card}}
\DeclareMathOperator{\Conv}{\mathrm{Conv}}

\crefname{subsection}{subsection}{the subsections}
\crefname{section}{section}{the sections}
\crefname{subsubsection}{subsubsection}{the subsubsections}

\crefname{Thm}{theorem}{theorems}
\Crefname{Thm}{Theorem}{Theorems}
\crefname{Prop}{proposition}{propositions}
\Crefname{Prop}{Proposition}{Propositions}
\crefname{Lemme}{lemma}{lemmas}
\Crefname{Lemme}{Lemma}{Lemmas}
\crefname{Cor}{corollary}{corollaries}
\Crefname{Cor}{Corollary}{Corollaries}
\crefname{Const}{construction}{constructions}
\Crefname{Const}{Construction}{Constructions}
\crefname{Ex}{exemple}{exemples}
\Crefname{Ex}{Exemple}{Exemples}
\crefname{Def}{definition}{definitions}
\Crefname{Def}{Definition}{Definitions}
\crefname{Not}{notation}{notations}
\Crefname{Not}{Notation}{Notations}
\crefname{Rem}{remark}{remarks}
\Crefname{Rem}{Remarque}{Remarques}
\usepackage{cancel}

\begin{document}
\title{Moduli spaces of quantum toric stacks and their compactification}
\author{Antoine BOIVIN}
\date\today
\maketitle

\begin{abstract}
    A toric variety is a normal complex variety which is completely described by combinatorial data, namely by a fan of strongly convex rational (with respect to a lattice) cones. Due to this rationality condition, toric varieties are (equivariantly) rigid since a deformation of the lattice can make it dense. 
A solution to this problem consists in considering quantum toric stacks. The latter is a stacky generalization of toric varieties where the "lattice" is replaced by a finitely generated subgroup of $\R^d$ (in the simplicial case as introduced by L.\textsc{Katzarkov}, E.\textsc{Lupercio}, L.\textsc{Meersseman} and A.\textsc{Verjovsky}).
 The goal of this paper is to explain the moduli spaces of quantum toric stacks and their compactification.
\end{abstract}

\newtheorem{Thm}{Theorem}[subsection]

\newtheorem{Prop}[Thm]{Proposition}

\newtheorem{Lemme}[Thm]{Lemma}

\newtheorem{Cor}[Thm]{Corollary}

\newtheorem{Const}[Thm]{Construction}

\newtheorem{DefProp}[Thm]{Definition-Proposition}

\theoremstyle{definition}

\newtheorem{Ex}[Thm]{Example}

\newtheorem{Def}[Thm]{Definition}

\newtheorem{Not}[Thm]{Notation}

\theoremstyle{remark}
\newtheorem{Rem}[Thm]{Remark}
\newtheorem{Avert}[Thm]{Warning}

\setcounter{secnumdepth}{4}
\setcounter{tocdepth}{3}
\tableofcontents

\section*{Introduction}

Quantum toric stacks are a (stacky) generalization of toric varieties introduced in \cite{katzarkov2020quantum} and \cite{boivin2020nonsimplicial}. They are analytic stacks built from the combinatorial data given by an "irrational" fan of cones, i.e. in any finitely generated subgroup of $\R^n$ (and no more in a lattice of $\R^n$ in the classical case). This construction is functorial and it induces an equivalence between the categories of fans and the one of quantum toric stacks. 

This permits to build moduli spaces of quantum toric stacks (it was not possible with toric varieties due to the rigidity of the rationality condition on the fan) by considering moduli spaces of fans.

The first goal of this paper is to precisely describe these moduli spaces, pursuing the work of \cite{katzarkov2020quantum}, when the combinatorics is fixed (in a sense made precise latter, in \cref{type_comb}). We prove in \cref{def} that they are orbifolds described by the quotient of a connected open semialgebraic subset by a finite group (the automorphism group of the combinatorial posets). The inequalities describing these semialgebraic subsets are given by signs of determinants coming from the combinatorics (see \cref{4-eq}). 
Moreover, in \cref{univ-family}, we prove the main feature of these moduli spaces: they are fine (in sense that each of them admits a universal family described thanks to the quantum GIT of \cite{katzarkov2020quantum}). In \cref{Ex}, we compute some families of examples of such moduli spaces: the moduli spaces containing the projective spaces of dimension $d$ is the quotient of an orthant of an affine space by an action of the permutation group $\Sfrak_{d+1}$ and the moduli spaces containing toric surfaces are quotients of a contractible space modded out by an action of a dihedral group.

The second goal of this paper is to describe a natural combinatorial compactification (in the sense that we can extend the universal family in a compatible way) of these moduli spaces: the boundary of such a compactification comes from degeneracies of the initial combinatorics. In \cref{5-compactf}, this compactification is realized thanks to an embedding of the previous semialgebraic subset in a Grassmannian manifold and by extending the action of the finite group to this Grassmannian manifold. Then, in \cref{extension_uf}, we extend the universal family by studying the degeneracies of the combinatorics occurring in the compactification. Finally, we provide explicit computations of all these degeneracies in the two previous families of examples (see \cref{5-EX}). In particular, the compactification of the moduli space of projective spaces of dimension $d$ arises as the quotient of a $d$-simplex in $\R\P^d$.

This work constitutes a first step towards the "toric big moduli conjecture" which states that the moduli spaces of quantum toric stacks with fixed combinatorics fit in a big moduli space of quantum toric stacks with fixed dimension and fixed number of generators which can be compactified (see \cite{boivin2022bigmoduli} for more details).

\section{Convention and notations}
\subsection{Quantum fans and quantum toric stacks}
In this subsection, we recall the needed defintions and theorems on quantum toric stacks (see \cite{katzarkov2020quantum} for the details of the constructions)
\begin{Def}
Let $\Gamma$ be a finitely generated subgroup of $\R^d$ such that $\Vect_\R(\Gamma)=\R^d$. A calibration of $\Gamma$ is given by:
\begin{itemize}
\item A group epimorphism $h \colon \Z^n \to \Gamma$
\item A subset $\Ical \subset \{1,\ldots,n\}$ such that $\Vect_\C(h(e_j), j \notin \Ical)=\C^d$ (this is the set of virtual generators)

\end{itemize}
This is a standard calibration if $\Z^d \subset \Gamma$, $h(e_i)=e_i$ for $i=1,\ldots,d$ and $\Ical$ is of the form $\{n-|\Ical|+1,\ldots,n\}$
\end{Def}

\begin{Def} \label{calib_q_fan_def}
A quantum fan $(\Delta,h : \Z^n \to \Gamma \subset \R^d,\Ical)$ in $\Gamma$ is the data of 
\begin{itemize}
\item a collection $\Delta$ of strongly convex polyhedral cones generated by elements of $\Gamma$ such that every intersection of cones of $\Delta$ is a cone of $\Delta$, every face of a cone of $\Delta$ is a cone and $\{0\}$ is a cone of $\Delta$.
\item a standard calibration $h$ with $\Ical$ its set of virtual generators
\item A set of generators $A$ i.e. a subset of $\{1,\ldots,N\} \setminus \Ical$ such that the 1-cone generated by the $h(e_i)$ for $i \in A$ are exactly the 1-cones of $\Delta$ 
\end{itemize}  
The fan is said simplicial if every cone of $\Delta$ is simplicial (i.e. which can be send on a cone of $\Cone(e_1,\ldots,e_k)$ by a linear automorphism of $\R^d$). \\
We note $\Delta(1)$ the cones of dimension 1 of $\Delta$ and $\Delta_{max}$ the maximal cones of $\Delta$. \\
With linear morphism which preserves inclusion of cones and calibration, the simplicial quantum fan forms a category denoted $\mathbf{QF}$.

\end{Def}

Every cone $\sigma=\Cone(h(e_{i_1}),\ldots,h(e_{i_{\dim(\sigma)}}))$ (isomorphic to $\Cone(e_1,\ldots,e_{\dim(\sigma)})$ by a linear morphism $L_\sigma$ and note $H$ a linear morphism described by a permutation $\chi$ of $\{1,\ldots,n\}$ such that $\chi(i_k)=k$ for $1 \leq k \leq \dim(\sigma)$) of a simplicial quantum fan $(\Delta,h,\Ical)$ describe an affine quantum toric stack (and, in particular, an analytic stack)
\[
\Uscr_\sigma \coloneqq [\C^{\dim(\sigma)}\times (\C^*)^{d-\dim(\sigma)}/\Z^{n-d}]
\]

where the action of $\Z^{n-d}$ on $\C^{\dim(\sigma)}\times (\C^*)^{d-\dim(\sigma)}$ is given by the following morphism 
\[
x \in \Z^{n-d} \mapsto  EL_\sigma h H^{-1}(0_{\R^{d}} \oplus x) \in (\C^*)^n
\]
where $E$ is the map \begin{equation} \label{def_exp}
    (z_1,\ldots,z_n) \in \C^n \mapsto \exp(2i\pi z_1,\ldots,\exp(2i\pi z_n)) \in (\C^*)^n
\end{equation}
In particular, with $\sigma=0$, we get the quantum torus $\Tscr_{h,\Ical}=[\C^*/\Z^{n-d}]$ (which is a Picard stack) which is dense in all these quantum toric stacks.

Thanks to the stabily by taking faces and stabily by intersection, a quantum fan $(\Delta,h,\Ical)$ define a diagram $D : \sigma \mapsto \Uscr_\sigma$ of affine quantum toric stacks. Indeed, for every $\sigma,\tau \in \Delta$, we have
\[
\begin{tikzcd}[ampersand replacement=\&,scale=0.05]
    \&\& {U_\sigma} \\
	{U_{\sigma\cap\tau}} \\
	\&\& {U_\tau}
	\arrow["{\text{open }}"{pos=0.5}, hook, from=2-1, to=1-3]
	\arrow["{\text{open}}"'{pos=0.5}, hook, from=2-1, to=3-3]
\end{tikzcd}
\]

\begin{Def}
    The quantum toric stack associated to the simplicial fan $(\Delta,h,\Ical)$ is the colimit $\Xscr_{\Delta,h,\Ical}$ of this diagram i.e. the gluing of the quantum toric stacks $\{\Uscr_\sigma\}_{\sigma \in \Delta}$ along their intersection. With stack morphisms which restrict on Picard stack morphism on torus, they form a category denoted $\mathbf{QTS}$.
\end{Def}

\begin{Thm}
    The correspondence $(\Delta,h,\Ical) \in \mathbf{QF} \mapsto \Xscr_{\Delta,h,\Ical} \in \mathbf{QTS}$ is an equivalence of categories.
\end{Thm}
This theorem permits us to study the moduli spaces of quantum toric stacks with the one of quantum fans.

\begin{Thm}[Quantum GIT]

If $(\Delta,h,\Ical)$ is a simplicial quantum fan then the quantum toric stack $\Xscr_{\Delta,h,\Ical}$ is a quotient stack
\[
\Xscr_{\Delta,h,\Ical}=[\Sscr(\Delta)/\C^{n-d}]
\]
where $\Sscr(\Delta)$ is a quasi-affine (classical) toric variety given by the combinatorics of $\Delta$ :
    \begin{equation} \label{Def_S}
        \Sscr(\Delta)=\bigcup_{\sigma \in \Delta_{max}} \C^{\sigma(1)} \times (\C^*)^{\sigma(1)^c} \subset \C^n;
    \end{equation}   
and $\C^{n-d}$ acts on $\Sscr$ by 
\[
t \cdot z=E(k \otimes id_\C(t))z
\]
where $k$ is a Gale transform of $h \otimes id_\R : \R^n \to \R^d$.
\end{Thm}

A Gale transform of $h \otimes id_\R$ is a morphism $k : \R^{n-d} \to \R^n$ such that 
\[\begin{tikzcd}[ampersand replacement=\&]
	0 \& {\R^{n-d}} \& {\R^n} \& {\R^d} \& 0
	\arrow[from=1-1, to=1-2]
	\arrow["k", from=1-2, to=1-3]
	\arrow["h", from=1-3, to=1-4]
	\arrow[from=1-4, to=1-5]
\end{tikzcd}\]
is exact

\subsection{Notations}
In these paper, we will use the following notations : 
\begin{itemize}
    \item $\sitename{Man}$ (resp. $\sitename{Top}$ )  are, respectively, the site of real manifolds (resp. compactly generated topological spaces) with Euclidean coverings; 
    \item $\catname{Poset}$ is the category of partially ordered sets with non-decreasing functions.
    \item $\Sfrak_I$ (resp. $\Sfrak_n$) is the group of bijections of $I$ (resp. permutations of $n$ objects) ; 
    \item $D_n$ is the dihedral group on $n$ objects ; 
    \item $P_\sigma$ is the linear isomorphism $\R^n \to \R^n$ associated to the permutation $\sigma \in \Sfrak_n$ i.e. $P=(\delta_{i,\sigma(j)})_{1 \leq i,j \leq n}$.
    
\end{itemize}

\section{Moduli spaces of quantum toric stacks}
 In this section, we study the moduli spaces of quantum toric stacks with fixed combinatorics. For this, we firstly introduce combinatoric type in order to give sense of the expression "fixed combinatorics". This permits us to define the desired moduli spaces which are global quotient stack of an open subspace $\Omega$ of a real affine space by the action of subgroup of a permutation group of the rays the combinatorics. These space $\Omega$ is in fact a connected semialgebraic subset whose equations are completely described by the starting combinatorics. Finally, we describe the universal families over these moduli spaces.

\subsection{Combinatorial type}
\label{type_comb}

\begin{Def} \label{def_typecomb_adm}
The combinatorial type of a fan $(\Delta, h : \Z^n \to \Gamma,\Ical)$ is the poset $\mathrm{comb}(\Delta)$ of subsets $I$ of $\{1,\ldots,n\}$ such that $\Cone(h(e_i), i \in I) \in \Delta$ ordered by inclusion. \\ Conversely, if $D$ is a poset, the calibration $(h,\Ical)$ is said $D$-admissible if there exists a quantum fan $(\Delta,h,\Ical)$ such that 
\[
\mathrm{comb}(\Delta) \simeq D
\]
\end{Def}

\begin{Ex}
The combinatorial type of the fan of the projective plane i.e. 
\[
\begin{cases}
\Delta=\{(e_1,e_2,-e_1-e_2,\Cone(e_1,e_2),\Cone(e_1,-e_1-e_2),\Cone(-e_1-e_2,e_2) \} \\
h \colon (x,y,z) \in \Z^3 \mapsto (x-z,y-z) \in \Z^2
\end{cases}
\]
is
\[
\mathrm{comb}(\Delta)=\{\{1\},\{2\},\{3\},\{1,2\},\{1,3\},\{2,3\}\}.
\]

\end{Ex}

Since a cone $\sigma=\Cone(\Acal) \subset \R^d$ is simplicial if, and only if,  for all $I \subset \Acal$, the cone $\Cone(I)$ is a face of $\sigma$ then we have the following statement:

\begin{Prop} \label{carac_comb_simp}
The property of being simplicial depends uniquely on the combinatorial type of the cone.
\end{Prop}

\begin{Lemme}
The correspondence $\mathrm{comb} : \catname{QF} \to {\Poset}$ is functorial.
\end{Lemme}

\begin{proof}
The poset morphism induced by a quantum fan morphism $(L,H) : (\Delta,h,\Ical) \to (\Delta',h',\Ical')$ is defined as follows : let $\sigma_I \in \Delta$. By definition of fan morphisms, $L\sigma$ is contained in a cone $\tau \in \Delta'$. The minimal cone for this property will be noted $\mathrm{comb}(L,H)(I)$. Hence, we have defined a poset morphism $\mathrm{comb}(L,H) : \mathrm{comb}(\Delta,h,\Ical) \to \mathrm{comb}(\Delta',h',\Ical')$ (i.e a non-decreasing map).
\end{proof}

\subsection{Definition and properties of moduli spaces}
\label{def}
\subsubsection{Definitions}

By definition of $D$-admissible calibration  (cf. \cref{def_typecomb_adm}), we get

\begin{Lemme} \label{lemme_comb}

If $h$ is $\comb(\Delta)$-admissible, the data of a calibrated fan $(\Delta,h,\Ical)$ is equivalent to the data of $(\comb(\Delta),h,\Ical)$ 
\end{Lemme}

We can defined the moduli space of quantum toric stacks with fixed combinatorics:

\begin{Def}
Let $D$ be the combinatorial type of a fan of $\R^d$, $n\geq\Card(D(1))$ and $\Ical \subset [\![1,n]\!] \setminus D(1)$. The moduli space of quantum toric stacks given by fans of combinatorial type $D$ and virtual generators $\Ical$ is the quotient stack\footnote{given by the stackification of the groupoid describing the equivalence relation} (over the site $\mathfrak{Man}_\R$) :
\begin{equation}
    \Mscr(d,n,D,\Ical) \coloneqq \{ h : \Z^{n} \to \R^d \mid  (h(e_i))_{i \in D(1)} \text{ is } D\text{-admissible} \}/\simeq \end{equation}
where $h \simeq h'$ if, and only if, there exists a fan isomorphism between the fan induced by $(D,h,\Ical)$ and the fan induced by $(D,h',\Ical)$ (cf \cref{lemme_comb}).
\end{Def}

\begin{Lemme}
Let $D$ be the combinatorial type of a fan of $\R^d$ which have a simplicial cone of maximal dimension $d$,  $n\geq\Card(D(1))$ and $\Ical \subset [\![1,n]\!] \setminus D(1)$. The stack $\Mscr(d,n,D,\Ical)$ is isomorphic to the stack
\[\{ h : \Z^{n} \to \R^d \mid h_{\mid \Z^d}=id_{\Z^d}, (h,\Ical) \text{ is } D\text{-admissible} \}/\simeq \]
or to the quotient of the subset $\Omega(D)$ of $\R^{d(n-d)}=(\R^d)^{n-d}$ (representing the different generators $h(e_i)$, $i \geq d$) mod out by isomorphisms of calibrated fans $(L,H)$ such that the set $I=\{j \mid \exists i \in \{1,\ldots,d\}, e_j=He_i \}\in D$ verifies \[\Cone(L(h(e_i)),i \in I)=\Cone(e_1,\ldots,e_d)
\]
\end{Lemme}

\begin{proof}
If $(\Delta,h,\Ical)$ has a simplicial cone $\sigma=\{v_1=h(e_{i_1},\ldots,v_d=h(e_{i_d})\}$ of maximal dimension then there exists a quantum fan isomorphisms $(L,H)$ which sends $(\Delta,h,\Ical)$ on a quantum fan which have a maximal cone which is $\Cone(e_1,\ldots,e_d)$. More explicitly, $L=(v_1 \ldots v_d)^{-1} : \R^d \to \R^d$ and $H$ is given by a permutation $\tau \in \Sfrak_n$ such that $\tau(i_k)=k$ for all $k \in \{1,\ldots,d\}$.
\end{proof}

\begin{Lemme} \label{Omega_ouvert}
If $D$ is simplicial then $\Omega(D)$ is a open subset of $\R^{d(n-d)}$.
\end{Lemme}

\begin{proof}
Since the fact that $d$ vectors of $\R^d$ are linearly independent is an open condition then the fact to be $D$-realizable for a simplicial combinatorial type $D$ is also an open condition.
\end{proof}

\begin{Rem}
The non-simplicial case is more complicated. Indeed, there exists combinatorial types of (non-simple) polytopes  which is not realizable over $\Q$ but over an algebraic extension of $\Q$ (see \cite{CP}). In these cases, the realization space is not open. We will study the polytopal case in \cref{4-polyt}.
\end{Rem}

\begin{Thm} \label{description_Mscr}
Let $D$ be the combinatorial type of a fan of 
$\R^d$ which have a simplicial cone of maximal dimension,  $n\geq\Card(D(1))$ and $\Ical \subset [\![1,n]\!] \setminus D(1)$. 
The moduli space of quantum toric stacks given by a fan of combinatorial type $D$ and with virtual generators $\Ical$ is isomorphic to the quotient stack
\begin{equation}
    \left[\Omega(D)/\Aut_{\catname{Poset}}(D) \times (\Z^n)^J \times \GL_J(\Z) \times \Sfrak_\Ical \right]
\end{equation}
where $J \coloneqq [\![1,n]\!] \setminus (D(1) \cup \Ical)$ is a set of generators which is neither appearing in the fan nor is virtual and $\Sfrak_\Ical$ is the group of bijections of $\Ical$. The action of $\Aut_{\mathbf{Poset}}(D) \times (\Z^n)^J \times \GL_J(\Z) \times \Sfrak_\Ical$ on $\Omega(D)$ is defined as follows :
Let $h \in\Omega(D)$, $\tau \in \Aut_\catname{Poset}(D), (\alpha_j)_{j \in J}$, $A \in \GL_J(\Z)$ and $\sigma \in \Sfrak_\Ical$. Let
\[
H=\begin{pmatrix}
    
    \begin{array}{c|c|c} 
      P_\tau &  (\alpha_i)_{i \in J}& 0 \\ 
      \hline 
      0& A &0 \\
      \hline
      0 & 0 & P_\sigma 
    \end{array} 
\end{pmatrix} \in \Mcal_n(\Z)
\]
and $L=(h(e_{\tau(1)}) \ldots h(e_{\tau(d)}))^{-1} $. Then,
\begin{equation}\label{action_esp_modules}
    (\tau,(\alpha_j)_{j \in J},A,\sigma) \cdot h=LhH^{-1} 
\end{equation}

\end{Thm}

\begin{proof}
A quantum fan isomorphism is the data of two linear morphisms $(L : \R^d \to \R^d, H : \R^n \to \R^n)$ which respects the cones and the virtual generators. Then, $H$ can be decomposed by blocks :
\begin{equation}
      H=\left(
    \begin{array}{c|c} 
      H_1 & 0\\ 
      \hline 
      0 & H_2
    \end{array} 
    \right)
\end{equation}

where $H_2 : \R^\Ical \to \R^{\Ical}$ is fully described by a permutation of $\Ical$. Now, it remains to us to describe $H_1$.
In the same way as \cite[section 11.2]{katzarkov2020quantum}, 
there exists a permutation $\tau$ such that $L=(h(e_{\tau(1)}) \ldots h(e_{\tau(d)}))^{-1} $ and $H_{\mid D(1)}=P_\tau$. More precisely, $\tau$ (and its inverse) must preserve the inclusions of cones. Hence it is a poset automorphism.
We deduce that $H_1$ is of the form
\[
H_1=\begin{pmatrix}
P_\tau & B \\
0 & A
\end{pmatrix}
\]
Since $H \in \GL_n(\Z)$, 
\[
\pm 1=\det(H)=\det(H_1)\det(H_2)=\det(P_\tau)\det(A)\det(P_\sigma)=\varepsilon(\sigma)\det(A)\varepsilon(\tau)
\]
Hence, $\det(A)=\pm 1$. We can conclude that $A \in\GL_J(\Z)$.
\end{proof}

\begin{Rem}
    Since an automorphism of posets preserves the 1-cones then the group $\Aut_\catname{Poset}(D)$ is a subgroup of $\Sfrak_n$.
\end{Rem}

\begin{Cor}
The stack $\Mscr(d,n,D,\Ical)$ is an orbifold if, and only if, the combinatorial type $D$ is maximal (i.e. $J=\emptyset$)
\end{Cor}

\begin{Not}
In what follows, we will suppose that the combinatorial type is maximal. Then we will omit the $\Ical$ in the notation (since $\Ical=[\![1,n]\!] \setminus D(1)$).
\end{Not}

\subsubsection{semialgebraicity}
\label{4-eq}
The goal of this paragraph is to prove that $\Omega(D)$ is a connected open semialgebraic subset of $\R^{d(n-d)}$ if $D$ is the combinatorial type of a simplicial fan.

\begin{Def}
Let $D$ be the combinatorial type of a complete simplicial fan of $\R^d$ with $n$ generators. The determinant associated to this combinatorial type is the map $\det_D : \R^{d(n-d)} \to \R^{D_{max}}$ defined by :
\[
\forall v=(v_{d+1},\ldots,v_{n}) \in (\R^{d})^{n-d}, {\det} _{D}(v)=(\det(v_i, i\in I))_{I \in D_{max}}
\]
where $v_1=e_1,\ldots,v_d=e_d$. \\
We will note $U(D)$ the (non-empty by hypothesis on $D$) open subset $\det_D^{-1}((\R^*)^{D_{max}})$ of $\R^{d(n-d)}$.
\end{Def}

In what follows, we will fix a combinatorial type $D$ which will be simplicial and all its maximal cones will be of cardinal $d$.

\begin{Rem} \label{inc_U_Om}
 .
Thanks to this hypothesis, $\Omega(D) \subset U(D)$
\end{Rem}

\begin{Lemme} \label{inegalite_omega}
The map $\varphi_D : U(D) \to \{ \pm 1\}^{\Delta_{max}}$ defined by
\[
\forall v \in U(D), \varphi_D(v)=\left(\frac{\det(v_i, i\in I)}{|\det(v_i, i\in I)|}\right)_{I \in D_{max}}
\]
is locally constant (i.e. constant on each connected component).
\end{Lemme}

\begin{proof}

It comes from the fact that the map $\varphi_D$ is continuous and the space $\{ \pm 1\}^{\Delta_{max}}$ is discrete.
\end{proof}

\begin{Prop} \label{ovuertfermeOmega}

The subset $\Omega(D)$ is a clopen subset in $U(D)$.
\end{Prop}

\begin{proof}

Since $D$ is simplicial then $\Omega(D)$ is an open subset of $\R^{d(n-d)}$ and hence by the inclusion of \cref{inc_U_Om}, an open subset of $U(D)$. Now we prove that it is also a closed subset of $U(D)$ : \\
Let $(h_n)$ be a sequence of $\Omega(D)$ which tends to $h$ in $\overline{\Omega}(D) \cap U(D)$. Since $h \in U(D)$ then, for $I \in D_{max}$ (which are of cardinal $d$)
\[
\det(h(e_i), i \in I) \neq 0
\]
For all $I \in D_{max}$, the cones $\sigma_I=\Cone(h(e_i), i \in I)$ are therefore simplicial cones. In the same way, the intersection properties are preserved thanks to the fact $h \in U(D)$. We deduce, thanks to \cref{carac_comb_simp}, that $h \in \Omega(D)$. It is conclude the proof of the closeness of $\Omega(D)$.
\end{proof}

Hence, the space $\Omega(D)$ is an union of connected component of $U(D)$. Therefore, thanks to \cref{inegalite_omega}, we get inequalities describing the connected components of $\Omega(D)$.

\begin{Cor} \label{ineq_omega}
Let $h_0 \in \Omega(D)$. The connected component of $\Omega(D)$ containing $h_0$ is given by the inequalities (in $h$) :
\[
\forall I \in D_{max}, \mathrm{sign}\det(h(e_i), i \in I)=\mathrm{sign}\det(h_0(e_i), i \in I)
\]
where $\mathrm{sign}$ is the map $\R \to \{\pm 1,0\}$ associated a real to his sign. The open $\Omega$ is a semialgebraic set of $\R^{d(n-d)}$.
\end{Cor}

\begin{Prop} \label{cnx_Omega}
If $D$ is the combinatorial type of a complete simplicial fan then $\Omega(D)$ is connected.
\end{Prop}

\begin{proof}
Since the sign of the determinant is fixed on the cone $(e_1,\ldots,e_d)$, the other are fixed inductively as well i.e. if $\sigma$ and $\tau$ are two maximal cones which have an ntersection of dimension $d-1$ then we can write $\sigma=\Cone(v_1,\ldots,v_{d-1},v_d)$ and $\tau=\Cone(v_1,\ldots,v_{d-1},v_{d+1})$ where
\[
\mathrm{sign}(\det(v_1,\ldots,v_{d-1},v_d))=-\mathrm{sign}(\det(v_1,\ldots,v_{d-1},v_{d+1}))
\]
\end{proof}

\subsubsection{Polytopal case}

\label{4-polyt}

In the literature, the space which is more frequently associated to a poset $D$ is the space describing the vertices of a polytope whose faces poset is isomorphic to $D$. 

\begin{Def}[\cite{richter1995realization}]
Let $P = \Conv(p_1,\ldots, p_n) \subset \R^d$ be a polytope of $\R^d$ with $n$ vertices and $\Bcal = (p_1, \ldots , p_{d+1})$ be a basis of the affine space $\R^d$.  The realization space $\Rcal(P, \Bcal)$ is the set of 
matrices $Q = (q_1, \ldots , q_n) \in \R^{dn}$ such that  \[\comb(\Conv(Q))=\comb(P)
\]
and
$q_i = p_i$ for $i = 1, \ldots , d + 1$.
\end{Def}
\begin{Rem}
The obtained space by change of basis are "equivalent" (in the sense of the end of the section)
\end{Rem}

\begin{Prop}[\cite{polytopal_fan}]
Let $D$ be the poset of faces of a (not necessarily simple) polytope $P \subset \R^d$ with $n$ vertices  and $\Bcal$ an affine basis of $\R^d$ (we will suppose that the first point of $\Bcal$ is $0$). Pose $\Pcal(P)$ thet set 
\[
\Pcal(P)\coloneqq \{(Q,v) \in \Rcal(P,\Bcal) \times \R^d \mid v \in \mathrm{Int}\Conv(Q)\}
\]
The map $(Q,v) \mapsto \{\text{vertices of }  Q-v\}$ is an open immersion $\Pcal(P) \hookrightarrow \Omega(D)$.

\end{Prop}

\begin{Rem}
The retired $v$ permits to the polytope to have 0 in the interior.
\end{Rem}

\begin{Avert}
The inclusion is always strict for $n>d+1$ even in dimension 2 i.e. even if every cones are generated by a face of a polytope, it does not imply that the chosen generators by the calibration are the vertices of a polytope  (see figure \ref{fig:Ce_faib_polytopaux}). 
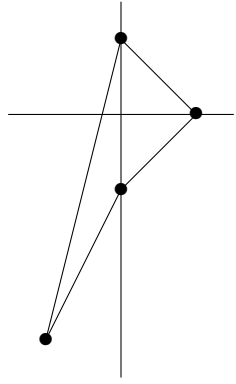
\begin{figure}[!ht]
    \centering
    \begin{tikzpicture}
    \draw (1,0) node{$\bullet$} -- (0,1) node{$\bullet$} -- (-1,-3) node{$\bullet$} -- (0,-1) node{$\bullet$} -- (1,0) ;
    \draw (0,0) -- (1.5,0);
    \draw (0,0) -- (0,1.5);
    \draw (0,0) -- (-1.5,0);
    \draw (0,0) -- (0,-3.5) ; 
    \end{tikzpicture}

    \caption{Counter-examples polytopal fan and calibrations}
    \label{fig:Ce_faib_polytopaux}
\end{figure}

\end{Avert}

The end of this subsection will be dedicated to prove that the spaces $\Omega(D)$ for $D$ non-simplicial can be as wild as possible.
\begin{Def}[stable equivalence]
\begin{itemize}
    \item Let $V \subset \R^n$ and $W \subset \R^{n+d}$ be two semialgebraic subsets such that $\pi(W)=V$ where $\pi$ is the projection $\R^{n+d} \to \R^n$. Then $V$ is a stable projection of $W$ if $W$ is of the form
    \[
    W=\{(v,v') \mid v \in V,  \forall i \in I,\forall j \in J,\varphi_i(v)\cdot v'>0, \psi_j(v) \cdot v'=0 \}
    \]
    where $I$ and $J$ are finite sets and the map $\varphi_i$ et $\psi_j$ are  polynomial maps $\R^n \to (\R^d)^*$
    \item Two semialgebraic subsets $V$ and $W$ are rationally equivalents if there exists a homeomorphism $f : V \to W$ such that $f$ and $f^{-1}$ are rational maps.
    \item The stable equivalence is the equivalence relation on the set of semialgebraic sets generated by the stable projection and the rational equivalence.
\end{itemize}
\end{Def}

\begin{Rem}

We can give a geometric interpretation to the stable projection : $V$ is a stable projection of $W$ if each fiber of the projection $W \to V$ is the interior of a polyhedron whose coefficients equations are given by polynomials on $V$.

\end{Rem}

\begin{Lemme}[\cite{richter2006realization}]
Let $V \subset \R^n$ and $W \subset\R^m$ two stably equivalent semialgebraic subsets.Then
\begin{itemize}
    \item If $k$ is a subfield of $\overline{\Q}$ and if $V \cap k^n=\emptyset$ then $W \cap k^m=\emptyset$
    \item $V$ and $W$ have "equivalent singularity structure". By example, $V$ is smooth if, and only if, $W$ is smooth.
\end{itemize}
\end{Lemme}

\begin{Thm}[of universality of \nom{Mnëv},\cite{Mnev1988},1986]
Every semialgebraic subset $V$ over $\Z$ is stably equivalent to the realization space of a polytope of dimension $d$ with $d+4$ edges.
\end{Thm}

\begin{Thm}[\nom{Richter}-\nom{Gebert},\cite{richter1995realization},1994]
Every semialgebraic subset $V$ over $\Z$ is stably equivalent to the realization space of a polytope of dimension $4$.
\end{Thm}

\begin{Cor}
The field $\overline{\Q} \cap \R$ is the smallest field where the 4-polytopes can be realized.
\end{Cor}

This results show that the sets $\Omega(D)$ can contain very singular open subsets and in order to realize some combinatorial type we need to order in an algebraic extension of $\Q$.

\begin{Rem}
An important consequent of these results is the fact that the polytopes realizing these combinatorial types are computationally hardly to handle : floats are just rational representation. Therefore, the approximation done by passing from reals to floats prevent the considered polytopes to have the good combinatorics.
\end{Rem}

\subsection{Examples}
\label{Ex}
\subsubsection{Projective spaces}
Let $n \in \N^*$ and $S_n$ the combinatorial type with $(n+1)$ rays and which maximal cones are the subfamilies with $n$ elements of $[\![1,n+1]\!]$ i.e. the combinatorial type of $\P^n$.

\begin{Lemme} \label{cnx_proj}
$\Omega(S_n)=\R_{<0}^n$ 
\end{Lemme}

\begin{proof}
We have the inclusion $\Omega(S_n) \subset \R_{<0}^n$ due to $(-1,\ldots,-1) \in \Omega(S_n)$ (since $\P^n$ has this combinatorial type). Moreover, $\Omega(D)$ is a connected subset of $ \R_{<0}^n$ (by \cref{cnx_Omega}) so $\Omega(S_n)=\R_{<0}^n$.

\end{proof}

\begin{Prop}\label{cohom_Mscr}
$\Mscr(n,n+1,S_n)$ has the homotopy type of the classifying space $B\Sfrak_{n+1}=K(\Sfrak_{n+1},1)$ and its cohomology groups are : 
\[
H^\bullet(\Mscr(n,n+1,S_n),\Z)=H^\bullet(\Sfrak_{n+1},\Z)
\]
\end{Prop}

\begin{proof}
    Since $\Omega(D) \to \Mscr(n,n+1,S_n)$ is a $\Sfrak_{n+1}$-fibration and $\Omega(S_n)$ is contractile then 
    \[
    \forall i \geq 0, \pi_i(\Mscr(n,n+1,S_n))=\pi_{i+1}(G)
    \]
    Since $G$ is discrete then 
    \[
     \pi_i(\Mscr(n,n+1,S_n))=\begin{cases}
         G \text{ if } i = 1 \\
         0 \text{ otherwise}
     \end{cases}
    \]
    i.e. $\Mscr(n,n+1,S_n)$ is homotopically equivalent to $K(G,1)=BG$ and hence \[
H^\bullet(\Mscr(n,n+1,S_n),\Z)=H^\bullet(B\Sfrak_{n+1},\Z)
\] Then, thanks to \cite[theorem 6.10.5]{weibel_1994}, we can deduce that it is isomorphic to the group cohomology ring $H^\bullet(\Sfrak_{n+1},\Z)$
\end{proof}

\begin{Ex}
For $n=2$, $\Sfrak_3=D_3$ and hence 
\[
H^k(\Mscr(n,n+1,S_n),\Z)=\begin{cases}
\Z \text{ if } k=0 \\
\Z/2\Z \text{ if }  k\equiv 2 [4] \\
\Z/6\Z \text{ if }  k\equiv 0 [4], k>0 \\
0 \text{ otherwise} 
\end{cases}
\]
\end{Ex}

\begin{Ex}[\cite{BS4}]
For $n=3$, $H^\bullet(\Mscr(3,4,S_3))$ is the graded ring described : 
\begin{align*}
    & \left\langle a, b,c,d \mid \deg(a) = 2, \deg(b) = 3, \deg(c) =3, \deg(d)=4, \right.\\
& \left. 2a = 2b = 4c = 3d = 0, ab^{2j} = a^{j+1}(c + a^2)^j \text{ for all } j \geq 1\right\rangle
\end{align*}
\end{Ex}

The computation of the group cohomology of the group $\Sfrak_n$ are not known for $n \geq 5$ but since $\Sfrak_n$ is finite then the rational cohomology group of $\Sfrak_n$ are easy to compute (see \cite[Chapter III : Proposition 10.1, Corollary 10.2]{book:29827961}:  

\begin{Prop}
For all $n \geq 1$, 
\[
H^k(\Mscr(n,n+1,S_n),\Q)=\begin{cases}
\Q \text{ if } k=0 \\
0 \text{ otherwise} 
\end{cases}
\]
\end{Prop}

\subsubsection{Dimension 2}
\begin{Lemme}
Let $D$ be the combinatorial type of a complete fan of $\R^2$. Then $D$ is isomorphic to
\[
(1,\ldots,n,(1,2),(2,3),\ldots,(|D(1)|-1,|D(1)|),(|D(1)|,1))
\]
\end{Lemme}
\begin{proof}

The isomorphism is given by ordonning the generators of the $1$-cones by their argument (by identifying $\R^2$ and $\C$).
\end{proof}

\begin{Not}
Note $C_n$ the combinatorial type \[
(1,\ldots,n,(1,2),(2,3),\ldots,(n-1,n),(n,1))
\]
which is the combinatorial type of a complete fan with $n$ generators in $\R^2$.
\end{Not}

\begin{Thm}\label{Mscr(2-)}
For all $n \geq 3$, the open subset $\Omega(C_n)$ is contractible.
\end{Thm}
    
\begin{proof}

By the inequalities of \cref{ineq_omega}, we have :
\begin{align*}
\Omega(D)=\{
v_3=(a_3,b_3),\ldots,v_{|D(1)|}=(a_{|D(1)|},b_{|D(1)|}) \in (\R^2 \setminus \{0\})^{|D(1)|-2}  \mid \\ a_3<0,b_{|D(1)|}<0,\det(v_i,v_{i+1})>0 \text{ pour } 3 \leq i \leq |D(1)|-1 \} 
\end{align*}
Consider now the image $\Omega'$ of $\Omega(D)$ by the map \[\pi \colon (x_1,\ldots,x_{[D(1)|-2]}) \in (\R^2 \setminus \{0\})^{|D(1)|-2} \mapsto \left(\frac{x_1}{|x_1|},\ldots,\frac{x_{D(1)|-2}}{|x_{D(1)|-2}|}\right) \in (\S^1)^{|D(1)|-2}.\] 
Since the cone are stable by the multiplication by a positive real number then $\pi_{|\Omega(D)} : \Omega(D) \to \Omega'$ is a $\R_{\geq 0}^{{|D(1)|-2}}$-fibration. We deduce that $\pi_{|\Omega(D)}$ is a homotopy equivalence.  Examine now the space $\Omega'$ : \\
By construction, the open subset $\Omega'$ is given by the same equations as $\Omega(D)$ : 
\begin{align*}
\Omega'=\{
v_3=(a_3,b_3),\ldots,v_{|D(1)|}=(a_{|D(1)|},b_{|D(1)|}) \in (\S^1)^{|D(1)|-2}  \mid \\ a_3<0,b_{|D(1)|}<0,\det(v_i,v_{i+1})>0 \text{ pour } 3 \leq i \leq |D(1)|-1 \} 
\end{align*}
This conditions on determinants become conditions on angles thanks to the following formula
\[
\det(u,v)=|\sin(\widehat{u,v})|
\]
for $u,v \in \S^1$ ($\widehat{u,v}$ is the angle between the two vectors). Moreover, we know that that quadrant $\R_{\geq 0}^2$ is already a cone of the fan. We deduce a homeomorphism
\[
\Omega'\simeq \left\{
(\alpha_2,\ldots,\alpha_{|D(1)|}) \mid 0 <\alpha_i<\pi/2 \text{ pour } 3 \leq i \leq |D(1)|-1, \sum_{i=2}^{|D(1)|} \alpha_i=3\pi/2
\right\}
\]
where $\alpha_2=\widehat{e_2,v_3},\ldots,\alpha_{|D(1)|}=\widehat{v_{|D(1)|},e_1}$. \\
Since $\alpha_{|D(1)|}$ depend on the other $\alpha_i$, we can project and hence we get a homeomorphism between $\Omega'$ and the space
\[
\left\{
(\alpha_2,\ldots,\alpha_{|D(1)|-1}) \mid 0 <\alpha_i<\pi/2 \text{ pour } 3 \leq i \leq |D(1)|-1,\pi/2<\sum_{i=2}^{|D(1)|-1} \alpha_i<3\pi/2
\right\}
\]
This space is clearly contractible and hence $\Omega(C_n)$ is contractible. 
\end{proof}

\begin{Prop}
    Let $n \geq 2$. The automorphism group of $C_n$ is the dihedral group $D_n$ on $n$ points.
\end{Prop}

\begin{proof}
Thanks to the normal fan construction, an automorphism of posets of $C_n$ is the same as an combinatorial automorphism of a $n$-gon i.e. an element of $D_n$. 
\end{proof}

\begin{Cor}[\cite{Handel} theorem 5.3]
Let $n \geq 3$ an odd integer. Then
\[
H^k(\Mscr(2,n,C_n),\Z)=H^k(D_n,\Z)=\begin{cases}
\Z &\text{ if } k=0 \\
\Z/n\Z \oplus (\Z/2\Z) &\text{ if } k \equiv 0[4], k \neq 0 \\
(\Z/2\Z) &\text{ if } k \equiv 2[4] \\
0 &\text{ otherwise} \\
\end{cases}
\]
Moreover, the graded ring $H^\bullet(\Mscr(2,n,C_n),\Z)$ is isomorphic to the polynomial ring
\[
\Z[a_2,d_4]/(2a_2,nd_4)
\]
where $a_2$ is a generator of $H^2(\Mscr(2,n,C_n),\Z)$ and $d_4$ is a  generator of $\Z/n\Z \subset H^4(\Mscr(2,n,C_n),\Z)$
\end{Cor}

\begin{Cor}[\cite{Handel} theorem 5.2]
Let $n \geq 3$ an even integer. Then 
\[
H^k(\Mscr(2,n,C_n),\Z)=H^k(D_n,\Z)=\begin{cases}
\Z &\text{ if } k=0 \\
\Z/m\Z \oplus (\Z/2\Z)^{\frac{k}{2}} &\text{ if } k \equiv 0[4], k \neq 0 \\
(\Z/2\Z)^{\frac{k-1}{2}} &\text{ if } k \equiv 1[4] \\
(\Z/2\Z)^{\frac{k+2}{2}} &\text{ if } k \equiv 2[4] \\
(\Z/2\Z)^{\frac{k-1}{2}} &\text{ if } k \equiv 3[4] \\
\end{cases}
\]
Moreover, the graded ring $H^\bullet(\Mscr(2,n,C_n),\Z)$ is isomorphic to the polynomial ring
\[
\Z[a_2,b_2,c_3,d_4]/\left(2a_2,2b_2,2c_3,nd_4,b_2^2+a_2b_2+\frac{n^2}{4}d_4,c_3^2+a_2d_4\right)
\]
where $a_2$, $b_2$ are two generators of $H^2(D_m,\Z)=\Z/2\Z \oplus \Z/2\Z$, $c_3$ is a generator $H^3(D_m,\Z)=\Z \/2\Z$ and $d_4$ is a generator of $\Z/m\Z \subset H^4(D_m,\Z)$
\end{Cor}

\subsubsection{\texorpdfstring{$\P^1 \times \P^2$}{P1xP2}}
Let $S_1=\{(1),(2)\}$ and $S_2=\{(3),(4),(5),(3,4),(4,5),(5,3)\}$ be the combinatorial type of, respectively, the projective line and the projective plane. Note $D$ the product combinatorial type  $D=D_1 \times D_2$. Hence, this is the combinatorial type of the fan of $\P^1 \times \P^2$.

\begin{Lemme}
$\Omega(D)$ is the set of $((a,b,c),(d,e,f)) \in (\R^3)^2$ such that
$\begin{cases}
a<0 \\
e<0 \\
f<0 \\
bd-ae<0 \\
dc-af<0
\end{cases}$
\end{Lemme} 
\begin{proof}

Same proof as \cref{cnx_proj}

\end{proof}
We have a projection $p :\Omega(D) \to \Omega(D_1) \times \Omega(D_2)=\R_{<0}^3$ defined by : 
\[
(a,b,c,d,e,f) \mapsto (a,e,f)
\]
Let $(a,e,f) \in \R_{<0}^3$. Then
\[
p^{-1}(a,e,f)\simeq \{(b,c,d) \in \R^3 \mid bd-ae<0, dc-af<0 \}
\]

\begin{Lemme}

For all $(a,e,f) \in \R_{<0}$, $\pi_{a,e,f} \colon (b,c,d) \in p^{-1}(a,e,f) \mapsto d \in \R$ is a homotopy equivalence of quasi-inverse $d \mapsto (0,0,d)$.
\end{Lemme}

\begin{proof}
Let $d \in \R$. Then 
\begin{align*}
   \pi_{a,e,f}^{-1}(d)&=\{(b,c) \in \R \mid bd-ae<0, dc-af<0 \} \\
   &=\begin{cases}
   ]- \infty,\frac{ae}{d}[\times]-\infty,\frac{ad}{e}[ \text{ si } d>0 \\
   \R^2 \text{ si } d=0 \\
    ]\frac{ae}{d},+\infty[\times]\frac{ad}{e},+\infty[ \text{ si } d<0 \\
   \end{cases}
\end{align*}
Since $ae>0$ then $(0,0)$ is in every fiber of $\pi_{a,e,f}$ which is therefore homotopic to $\{(0,0)\}$.

\end{proof}

We deduce from this

\begin{Prop}
The projection $\Omega(D) \to \R_{<0}^3 \times \R$ defined by :
\[
(a,b,c,d,e,f) \mapsto (a,d,e,f)
\]
is a homotopy equivalence of quasi-inverse $(a,d,e,f) \mapsto (a,0,0,d,e,f)$.
\end{Prop}

\begin{Cor}
The open subset $\Omega(D)$ is contractile and $\Mscr(3,5,D)$ has the homotopy type of $B\Aut_{\Poset}(D)=B(\Z/2\Z \times \Sfrak_3)=BD_6$.
\end{Cor}

In the same manner as \cref{cohom_Mscr}, we have : 
\begin{Cor}
We have the following equalities
\[
H^k(\Mscr(3,5,D),\Z)=\begin{cases}
\Z &\text{ if } k=0 \\
\Z/2\Z &\text{ if }  k\equiv 2 [4] \\
\Z/12\Z &\text{ if }  k\equiv 0 [4], k>0 \\
0 &\text{ otherwise}

\end{cases}
\]
\end{Cor}

\subsection{Universal family}
\label{univ-family}

In this subsection, we will describe an universal family of quantum toric stacks associated to $\Mscr$. The first step is the embedding of $\Omega$ in a  Grassmannian manifold. This embedding have two consequences : the first one is the simplification of the action \eqref{action_esp_modules} of the group $G$ on the image of $\Omega$ and the second one is this induces a natural compactification since the Grassmannian is compact. \\

We will begin some recall on Grassmannian manifolds:

\begin{Prop}
The space $\Gr(n-d,\R^n)$ of subvector spaces of dimension $n-d$ of $\R^n$ verifies the following statements
\begin{itemize}
    \item It is a compact differentiable manifold of dimension $d(n-d)$ whose charts are given by the open subsets $\{U_I\}_{|I|=d} $ of vector spaces transverse to $\R^I \oplus 0 \subset \R^n$ which are isomorphic to $\R^{d(n-d)}$ (we will note by $\varphi_{I} : U_I \to \R^{d(n-d)}$ this isomorphism) ;  \\
    \item The kernel map $\ker : \mathrm{Epi}(\R^n,\R^d) \to \Gr(n-d,\R^n)$, between the space of linear epimorphisms and the Grassmannian manifold, is a $\GL_d(\R)$-fiber bundle. Dually, 
    the image map $\ker : \mathrm{Mono}(\R^{n-d},\R^n) \to \Gr(n-d,\R^n)$, between the space of linear monomorphisms and the Grassmannian manifold, is a $\GL_{n-d}(\R)$-fiber bundle. \\
    Explicitly, for $I=\{i_1<\ldots<i_d\}$ a set of cardinal $d$ (note $I^c=\{i_1^c < \ldots <i^c_{n-d}\}$ its complementary), if we note $\mathrm{Epi}_I(\R^n,\R^d)$ the open subset of $\mathrm{Epi}(\R^n,\R^d)$ with epimorphism which restricts to an isomorphism $\R^I \to \R^d$, the map $\ker$ restricts to a map $\mathrm{Epi}_I(\R^n,\R^d) \to U_I \simeq (\R^d)^{n-d}$ which is a trivial fiber bundle thanks to the section $s_{\mathrm{Epi},I}$ defined by
    \begin{equation} \label{section_epi}
        \forall (v_1,\ldots,v_{n-d}) \in (\R^d)^{n-d}, s_{\mathrm{Epi},I}(v_1,\ldots,v_{n-d}) : x \in \R^n \mapsto \sum_{k=1}^d x_{i_k}e_{i_k}+\sum_{j=1}^{n-d} y_{i^c_j} v_{i^c_j}. 
    \end{equation}
    In the same way, for $I=\{i_1<\ldots<i_d\}$ a set of cardinal $d$, if we note $\mathrm{Mono}_{I^c}(\R^{n-d},\R^n)$ the open subset of $\mathrm{Mono}(\R^{n-d},\R^n)$ containing monomorphisms $k=(k_1,\ldots,k_n)$ which corestricts to an isomorphism \[(k_{i^c_1},\ldots,k_{i^c_{n-d}})\colon\R^{n-d}\to \R^{I^c},\]
    the map $\im$ restricts to a map $\mathrm{Mono}_{I^c}(\R^{n-d},\R^n) \to U_I \simeq (\R^{n-d})^d$ which is a a fiber trivial bundle thanks to the section $s_{\mathrm{Mono},I}$ defined by 
    \begin{equation}
    \forall (w_1,\ldots,w_{d}) \in (\R^{n-d})^d, s_{\mathrm{Mono},I}(w_1,\ldots,w_d) : x \mapsto \R^{n-d} \mapsto \sum_{k=1}^{n-d} x_{i^c_k} e_{i_k} + \sum_{l=1}^d \left\langle x,v_i\right\rangle e_{i_l}
    \end{equation}
    \item The following diagram commutes
    \begin{equation} \label{iso_grassm}
        \begin{tikzcd}[ampersand replacement=\&]
	{\mathrm{Epi}_I(\R^n,\R^d)/\GL_d(\R)} \&\& {U_I \subset \Gr(n-d,\R^n)} \\
	\\
	\&\& {\mathrm{Mono}_I(\R^{n-d},\R^n)/\GL_{n-d}(\R)}
	\arrow["\ker", from=1-1, to=1-3]
	\arrow["\im"', from=3-3, to=1-3]
	\arrow["{\text{Gale transform}}"', from=1-1, to=3-3]
	\arrow["\simeq"', from=1-1, to=1-3]
	\arrow["\simeq", from=3-3, to=1-3]
	\arrow["\simeq", from=1-1, to=3-3]
\end{tikzcd}
    \end{equation}
    
\end{itemize}
\end{Prop}

\begin{Lemme} \label{equivariance}
The morphism $\varphi_{[\![1,d]\!]}^{-1} : \R^{d(n-d)}\to \Gr(n-d,\R^n)$ is an $G$-equivariant open immersion for the action \eqref{action_esp_modules} at left and the action
defined by, for $V^{n-d}$ a subvector space of $\R^n$ and $(\sigma,\tau) \in G$, 
\begin{equation} \label{action_grass}
     (\sigma,\tau) \cdot V=\begin{pmatrix} P_\sigma^{-1} & 0 \\ 0 & P_\tau^{-1}
 \end{pmatrix} V
\end{equation}

\end{Lemme}

\begin{proof}
Up to isomorphisms (given by \eqref{iso_grassm}), it is the inclusion \[\mathrm{Epi}_I/\GL_{d}(\R)\hookrightarrow \mathrm{Epi}/\GL_d(\R).\] 
A quantum fan isomorphism is given by the following commutative diagram
\begin{equation} \label{diagram_morphism_quantique}
    \begin{tikzcd}
	{\R^n} & {\R^n } \\
	{\R^d} & {\R^d}
	\arrow["h"', from=1-1, to=2-1]
	\arrow["L"', from=2-1, to=2-2]
	\arrow["H", from=1-1, to=1-2]
	\arrow["{h'}", from=1-2, to=2-2]
\end{tikzcd}
\end{equation}
where $L$ and $H$ are linear isomorphisms (verifying some properties on the quasi-lattices).

Since the morphism $H$ sends the cones of $\Delta_h$ on cones of $\Delta_{H'}$, it permutes generators of the cones. So we can write it as \eqref{action_grass} for a well-chosen pair $(\sigma,\tau)$.

The morphism $L$ is fully determined by the permutation $\sigma$ and is equal to the matrix $L_\sigma$ of the action \eqref{action_esp_modules}. 
By the diagram \eqref{diagram_morphism_quantique}, we get the equality
\[
\ker(h')=H(\ker(h))
\]
(since $L$ is an isomorphism). In other words, we have
\[
\ker((\sigma,\tau) \cdot h)=(\sigma,\tau)\cdot \ker(h)
\]
where the left action is the action \eqref{action_esp_modules}. In other words, the morphism $\iota$ is $G$-equivariant.
\end{proof}

We will note again $\Omega$ the image of $\Omega$ by $\varphi^{-1}_{[\![1,d]\!]}$. Hence, we can interpret $\Mscr$ as a open substack of $[\mathrm{Gr}_\R(n-d,\R^n)/G]$.

\begin{Not}
In what follows, we will note $[k]$ an element of the Grassmannian manifold seen as an equivalence class of morphisms $\R^{n-d} \to \R^n$. 
\end{Not}

By lifting this action with the section $s_{\mathrm{Mono},I}$, we get the following statement :

\begin{Lemme}
The action of $G$ on the Grassmannian manifold comes from the descent to quotient of the actions of $G$ on $\mathrm{Mono}(\R^{n-d},\R^n)$, defined, respectively, by
\[
(\sigma,\tau) \cdot k=\begin{pmatrix} P_\sigma^{-1} & 0 \\ 0 & P_\tau^{-1} \end{pmatrix} k
 \]
\end{Lemme}

\begin{proof}
This comes from the diagram \eqref{diagram_morphism_quantique} and the fact that we will consider the kernel of $h$ and the image of $k$.
\end{proof}

The goal of the end of this section is to prove the following theorem

\begin{Thm} \label{famille_universelle_annonce}
    The moduli stack $\Mscr$ admits a universal family i.e. there exists a stack morphism $\Xscr \to \Mscr$ whose fibers are quantum toric stacks of combinatorial type $D$. 
\end{Thm}

\begin{Not}
For every $I$ of cardinal $n-d$ and every element $[k] \in \Omega' \cap U_I \eqqcolon \Omega_I$, we can define an action of $\C^{n-d}$ on $\Sscr$ : 
\[
t \cdot_k z\coloneqq E(s_{\mathrm{Mono},I}([k])(t))z
\]
\end{Not}
We can make this construction in family by considering the quotient
\begin{equation}
\Xscr_I=[
\Sscr \times \Omega_I/\C^{n-d}
]
\end{equation} 
where the action of $\C^{n-d}$ on $\Sscr \times \Omega_I$ is defined by
\[
t \cdot (z,k)=(E(s_{\mathrm{Mono},I}([k])(t))z,k)
\]  
In particular, $\Xscr=\Xscr_{[\![d+1,n]\!]}=[\Sscr \times \Omega/\C^{n-d}]$.

The projection on the second coordinate $\Sscr \times \Omega_I \to \Omega_I$ (which is $\C^{n-d}$-invariant) descend to quotient into a morphism
\[
\Xscr_I \to \Omega_I
\]

\begin{Lemme}

Let $I$ be a set of the form $[\![1,n]\!] \setminus \sigma \cdot [\![1,d]\!]$ where $\sigma \in G$. The action of the group $G$ on the Grassmannian manifold $ \Gr(n-d,\R^n)$ which restricts on an action on $\Omega_I$.
\end{Lemme}

\begin{proof}
Let $\sigma \in G$ and $[h] \in \Omega_I$. The open subset $\Omega_I$ is stable by $G$ if, and only if, for all $\tau \in G$, the matrix of $h \circ \tau_{|\R^{I^c}}$ is invertible. This is  equivalent to ask that the matrix of $h \circ \tau\sigma_{|\R^{d}}$ is invertible, which is always true by definition of $G$ (since the action of $\tau\sigma \in G$ sends the cone $\Cone(e_1,\ldots,e_d)$ on a simplicial cone of dimension $d$).
\end{proof}

\begin{Rem}
This is not true for a general $I$ : \\
Let $h : \Z^4 \to \R^2$ a morphism defined by :
\[
h(e_1)=e_1, \ h(e_2)=e_2, \ h(e_3)=-e_1, \ h(e_4) \in \R_{\neq 0} \times \R_{<0}
\]
Then $[h] \in \Omega_{\{2,4\}}$ and $(1 \ 2 \ 3 \ 4 ) \cdot [h] \notin \Omega_{\{2,4\}}$.

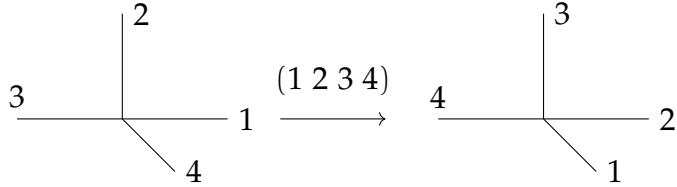
    \begin{tikzpicture}[scale=0.7]
    \draw (0,0) to (0,2) node[right]{$2$} ;
    \draw (0,0) to (2,0) node[right]{$1$} ;
    \draw (0,0) to (-2,0)node[right, above]{$3$}   ;
    \draw (0,0) to (1,-1)node[right]{$4$}   ;
    \draw[->] (3, 0) to (5,0)
    node[above] at (4,0.2) {$(1\ 2\ 3 \ 4)$} ; 
    \draw (8,0) to (8,2) node[right]{$3$} ;
    \draw (8,0) to (10,0) node[right]{$2$} ;
    \draw (8,0) to (6,0)node[right, above]{$4$}   ;
    \draw (8,0) to (9,-1)node[right]{$1$}   ;
\end{tikzpicture}
\end{Rem}
Since the end of this subsection, we will suppose $I=\{i_1<\ldots<i_{n-d}\}$ is of the form $[\![1,n]\!] \setminus \sigma \cdot [\![1,d]\!]$. \\
This action extend on $\Omega_I \times \Sscr$ by :
\begin{equation} \label{action_locale}
\sigma \cdot (k,z)=(\sigma \cdot [k],(z_{\sigma(1)},\ldots,z_{\sigma(n)}))
\end{equation}
since $\sigma \in G$. \\

In what follows, we will study the behavior of this action with respect to the action of $\C^{n-d}$ on $\Omega_I \times \Sscr$ in order to make this action descends to the quotient.

Firsty, we will see the relation with the section $s_{\mathrm{Mono},I}$ : 

\begin{Lemme} \label{prep_equiv}
For all $[k] \in \Omega_I$ and all $\sigma \in G$, we have :
\begin{equation}\label{def_K_sigma}
  s_{\mathrm{Mono},I}(\sigma \cdot [k])=(\sigma \cdot (s_{\mathrm{Mono},I} ([k]))) \circ K_{\sigma,I}^{-1}([k])  
\end{equation}
where $K_{\sigma,I}([k])=((\sigma_{\mathrm{Epi},I} ([k]))_{i_m}(e_{\sigma(j)}))_{1 \leq j,m \leq n-d} \in \GL_{n-d}(\R)$
\end{Lemme}

\begin{proof}
Let $[k] \in \Omega_I$, $k=s_{\mathrm{Mono},I}([k])$ and $\sigma \in G$.
By definition, the morphism $\sigma \cdot k=(k_{\sigma(1)},\ldots,k_{\sigma(n)})$ is in the orbit of $\sigma \cdot k \circ K_{\sigma,I}([k])^{-1}$ which is in the image of $s_{\mathrm{Mono},I}$ (by construction of $K_{\sigma,I}([k])$). Consequently,
\[
s_{\mathrm{Mono},I}(\sigma\cdot [k])=\sigma \cdot k \circ K_{\sigma,I}([k])^{-1}=\sigma \cdot s_{\mathrm{Mono},I}([k]) \circ K_{\sigma,I}([k])^{-1}
\]

\end{proof} 

\begin{Prop} \label{etape_calcul}
For all $[k] \in \Omega_I$, $z \in \Sscr$, $t \in \C^{n-d}$ and $\sigma \in G$, 
\[
\sigma \cdot (k,t \cdot z)=(K_{\sigma,I}([k]) t) \cdot (\sigma\cdot (k,z) ) 
\]
\end{Prop}

\begin{proof}
Let $[k]\in \Omega_I$, $z \in \Sscr$, $t \in \C^{n-d}$ et $\sigma \in G$. Then, if we note $z_\sigma$ the point $(z_{\sigma_1},\ldots,z_{\sigma(n)})$, we have, thanks to \cref{prep_equiv} : 
\begin{align*}
(K_{\sigma,I}([k]) t) \cdot (\sigma\cdot([k],z)  ) &=(K_{\sigma,I}([k]) t) \cdot ([k_\sigma],z_\sigma) \\
&=(k_\sigma,E(s_{\mathrm{Mono},I}([k_\sigma])(K_{\sigma,I}([k])t))z_\sigma) \\
&=\sigma \cdot (k,E(s_{\mathrm{Mono},I} ([k])(t))z)=\sigma \cdot (t \cdot (k,z))
\end{align*} 
\end{proof}
The matrices $K_{\sigma,I}$ verify a cocycle collection of the form \eqref{eq_2-comm} :

\begin{Prop}

Let $I$ be a set of cardinal $n-d$, $\sigma, \tau \in G$ et $[k] \in \Omega_I$.
\begin{equation} \label{eq_2-comm_Aut}
 K_{\sigma,I}(\tau \cdot [k])(K_{\tau,I}([k]))=K_{\sigma\tau,I}([k])
\end{equation}
\end{Prop}

\begin{proof} 
By \eqref{def_K_sigma}, for all $[k] \in \Omega_I$ and all $\sigma,\tau \in G$, we have :
\begin{align*}
    (\sigma \tau s_{\mathrm{Mono},I}([k]))K_{\sigma\tau,I}^{-1}&=s_{\mathrm{Mono},I}(\sigma \tau \cdot [k]) \\
    &=\sigma \cdot s_{\mathrm{Mono},I}(\tau \cdot [k]) K_{\sigma,I} (\tau \cdot [k])^{-1} \\
    &=\sigma\tau \cdot s_{\mathrm{Mono},I}([k]) K_{\tau,I}([k])^{-1} K_{\sigma,I} (\tau \cdot [k])^{-1} 
\end{align*}
By injectivity of $k$, we get
\[
K_{\sigma\tau,I}^{-1}=K_{\tau,I}([k])^{-1} K_{\sigma,I} (\tau \cdot [k])^{-1}
\]
By passing to the inverse the previous equality, we get the desired equality.
\end{proof}

We deduce from this two results, the \cref{descente_quotient} and the \cref{descente_action}the following result:  
\begin{Prop} \label{pass_quotient}

The action of $\sigma$ on $\Omega_I$ descend to quotient on an action on $\Xscr_I$.
\end{Prop}
Moreover, the projection $\Xscr_I \to \Omega_I$ is $G$-equivariant. This descends to quotient:
\begin{equation} \label{proj_famille}
\Xscr_I/G \to [\Omega_I/G]
\end{equation} 
In conclusion, we get the universal family of \cref{famille_universelle_annonce}
\begin{Thm}
The projection 
\begin{equation} \label{famille_universelle}
\Xscr/G \to [\Omega_{[\![d+1,n]\!]}/G] \simeq \Mscr
\end{equation}
 is the universal family of quantum toric stacks of combinatoric type $D$.
\end{Thm}

\section{Compactification}

\subsection{Compactification of moduli spaces}
\label{5-compactf}

In this subsection,we will construct a natural compactification of $\Mscr$ i.e. describe a stack $\overline{\Mscr}$ which contains $\Mscr$ as dense open substack with a family $\overline{\Xscr} \to \overline{\Mscr}$ which extends the family \eqref{famille_universelle}.

\begin{Def}
Let $\Xscr$ be a topological stack. Let $\Ascr$ be a substack of $\Xscr$. The closure $\overline{\Ascr}$ of $\Ascr$ in $\Xscr$ is the closed substack satisfying the following universal property: \\
For all closed substack $\Fscr$ of $\Xscr$ containing $\Ascr$, there exists a monomorphism $\overline{\Ascr} \hookrightarrow \Fscr$ making 2-commutes the following diagram: 
\[\begin{tikzcd}
	& {\overline{\Ascr}} \\
	\Ascr & {} && \Xscr \\
	& \Fscr
	\arrow[hook, from=1-2, to=2-4]
	\arrow[hook, from=2-1, to=1-2]
	\arrow[hook, from=2-1, to=3-2]
	\arrow[hook, from=3-2, to=2-4]
	\arrow[hook', from=1-2, to=3-2]
\end{tikzcd}\]
\end{Def}

Thanks to the correspondance between the closed substack of $[X/G]$ and the $G$-invariant closed subset of $X$, we deduce 

\begin{Prop}
Let $A$ be a $G$-invariant subset of $X$ and note $\overline{A}$ the $G$-invariant closure of $A$ in $X$ i.e. the smallest $G$-invariant closed subset containing $A$. Then the closure of $[A/G]$ is $[\overline{A}/G]$.
\end{Prop}

When $G$ is finite (this will be the case in this paper), we have an explicit description dof $G$-invariante closure : 

\begin{Lemme}
Let $G$ be a finite group, $X$ a topological space with a continuous action of $G$ and $A$ a $G$-invariant subset of $X$ . Then the $G$-invariante closure of $A$ in $X$ is the closed subset 
\[
\bigcup_{g \in G} g \cdot \mathrm{Cl}(A)
\]
where $\mathrm{Cl}(A)$ is the closure of $A$ in $X$.
\end{Lemme}

Suppose $\Xscr=[X/G]$ is a quotient differentiable stack and $\Ascr=[A/G]$ is an open substack $\Xscr$. Then the closure $\overline{\Ascr}=[\overline{A}/G]$ of $\Ascr$ in $\Xscr$ (seen as stacks over $\sitename{Top}$) can be seen as stack over $\sitename{Man}$ if we consider the closure $\overline{A} \stackrel{j}{\hookrightarrow} X$ as the sheaf $j^*\underline{X}$ over $\sitename{Man}$ i.e.
\[
j^*\underline{X}(U)=\left\{ f : U \to X \mid f(U)\subset\overline{A} \right\}
\]
and $\overline{\mathscr{X}}=[\overline{A}/G]$ as the quotient stack $j^*\overline{X}/G$ (cf. \cite{Romagny}).

\begin{Def}
We call topological closure of the differentiable stack $\Ascr=[A/G]$ in the differentiable stack $\Xscr=[X/G]$ the stack $j^*\overline{X}/G$ defined over $\sitename{Man}$.
\end{Def}

Let $\overline{\Omega}$ be the $G$-invariant closure of $\Omega$ in $\mathrm{Gr}_\R(n-d,\R^n)$.
\begin{Thm}
The space $\overline{\Omega}$ is the connected closed semialgebraic subset of $\mathrm{Gr}(n-d,\R^n)$ (seen as  real submanifold of $\R^{n^2}$, cf. \cite[Théorème 3.4.4]{bochnak2013real}).
\end{Thm}

\begin{proof}
Since the closure of semialgebraic set is semialgebraic (see \cite[Proposition 2.2.2]{bochnak2013real}), $\overline{\Omega} \cap U_i$ is a semialgebraic subset of $U_i \simeq \R^{d(n-d)} \subset \Gr(n-d,\R^n)$ which can be (algebraicly) embedded in $\R^{n^2}$. By taking the (finite) union on the $I$, we deduce 
\[
\overline{\Omega}=\bigcup_I \overline{\Omega}\cap U_i
\]
is a semialgebraic subset. \\
The $\overline{\Omega}$ is connected thanks to the connectedness of $\Omega$ (since we are in the simplicial case, see \cref{cnx_Omega}).
\end{proof}

More precisely, we have the following result: 

\begin{Prop} \label{cnx_eq_ombar}
Let $\varepsilon=\det_D(h) \in \{\pm 1\}^{\Delta_{max}} $ with $h \in \Omega$ (this do not depend of the choice of $h$ in $\Omega$). Then we have 
\[
\Omega=\{
h=(h_{d+1},\ldots,h_n) \in \R^{d(n-d)} \mid \forall I \in D_{max}, 
\mathrm{sign}\det(h_{d+1},\ldots,h_n)=\varepsilon_I\}
\]
The closed subset $\overline{\Omega}$ (embedded in $\P\left(\bigwedge^r \R^n\right)$ by the Plücker embedding) is a connected component of 
\begin{equation} \label{eq_ombar}
    \left\{[p_I] \in \Gr(n-d,\R^n) \hookrightarrow \P\left(\bigwedge^r \R^n\right) \mid \forall I,J \in D_{max}, \mathrm{sign}(p_Ip_J)\in\overline{\varepsilon_I}\overline{\varepsilon_J}  \right\}
\end{equation}
with equality if \eqref{eq_ombar} is connected (here, $\mathrm{sign}(x)  \in \overline{1}$ if $x \geq 0$ and $\mathrm{sign}(x)  \in \overline{-1}$ if $x \leq 0$).

\end{Prop}
\begin{Rem}
The sign which appear here are well-defined because the change of representative $\lambda p \simeq p$ multiply the numbers by $\lambda^2$ and hence do not change the sign. 
\end{Rem}

\begin{proof}
Note $F$ the right hand side. It is a closed subset of $\P\left(\bigwedge^r \R^n\right)$. Moreover, we have an injection
\[
\Omega \hookrightarrow F
\]
since \{1,\ldots,d\} is a maximal cone of $D$ and the associated determinant $p_{1,\ldots,d}$ is $1 \geq 0$. The other equations are immediately verified.
\end{proof}

\subsection{Extension of the universal family}
\label{extension_uf}
The goal of this subsection is to extend the universal family \eqref{famille_universelle} into a family over $[\overline{\Omega}/G]$.\\
Note $\overline{\Omega}_I$ the intersection $\overline{\Omega} \cap U_I$ which is the closure of $\Omega_I=\Omega \cap U_I$ in $U_I \simeq \R^{d(n-d)}$. The closed subspace $\overline{\Omega_I} \setminus \Omega_I$ correspond to the configurations of vectors which is no more $D$-admissible because some cones of $D$ is no more strongly convex. For solving this problem, we will consider the degenerated combinatorial types from $D$ i.e. combinatorial type where we have removed some cones which is not more strongly convex for a calibration $h \in \Omega(D)$. More precisely,

\begin{Def}
Let $D$ be the combinatorial type of a fan in $\R^d$. A degenerated combinatorial type of $D$ is a sub-poset $D'$ of $D$ which have the same 1-cones as $D$, which is stable by intersection (i.e. if $\sigma, \tau \in D'$ then $\sigma \cap \tau \in D'$),  which satisfies 
\[
\forall \sigma \in D', \forall \tau \in D,  \tau\leq \sigma \Rightarrow \tau \in D'
\]
(to abbreviate, we will say the poset $D'$ is stable by taking the faces, by analogy with the faces of a cone) and dont the maxima of the cardinals of the elements of these two combinatorial types are equal.
\end{Def}

\begin{Ex}
Let $D$ be the combinatorial type of the fan of the projective plane i.e.
\[
D=(1,2,3,\{1,2\},\{2,3\},\{3,1\})
\]
Then the degenerated combinatorial type from $D$ are the following posets:
\begin{align*}
 &D_1'=(1,2,3,\{1,2\},\{3,1\}) \\
&D_2'=(1,2,3,\{1,2\},\{2,3\})\\
&D_3'=(1,2,3,\{2,3\},\{3,1\})\\
&D_{1,2}'=(1,2,3,\{1,2\})\\
&D_{2,3}'=(1,2,3,\{2,3\})\\
&D_{3,1}'=(1,2,3,\{3,1\})   
\end{align*}

These degenerated combinatorial type correspond to subfamilies of the cones of $\Delta$ where some cones are not strongly convex and hence we have to remove them in order to have a fan. By example, for the combinatorial type $D_1'$, we get fans of the form \eqref{fig:fandegen1}.
 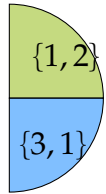
\begin{figure}[!ht]
    \centering
    \begin{tikzpicture}[scale=0.5]
    \tikzmath{\ray = 2.5;}
    \filldraw[fill opacity=0.5,fill=couleur4] (0,0) -- (\ray,0) arc (0:90:\ray cm) -- cycle;
    \filldraw[fill opacity=0.5,fill=couleur3] (0,0) -- (0,-\ray) arc (270:360:\ray cm) -- cycle;
    \node at (.35,.35) [above right] {$\{1,2\}$};
    \node at (0,-2) [above right] {$\{3,1\}$};
\end{tikzpicture}
    \caption{\'Eventail associé à $D'_1$}
    \label{fig:fandegen1}
\end{figure}
\end{Ex}

\begin{Lemme}
A degenerated combinatorial type of a simplicial combinatorial type is simplicial.
\end{Lemme}

\begin{proof}
It follows from \cref{carac_comb_simp}.
\end{proof}

By the definition of $\Sscr$ given in \eqref{Def_S}, we get : 

\begin{Lemme}
If $D'$ is a degenerated combinatorial type of $D$ then
\begin{align}
\label{var_deg_0}
\Sscr(D')&=\Sscr(D) \setminus \bigcup_{\sigma \in D \setminus D'} (\C^*)^{\sigma^c} \times 0^{\sigma} \\
&=\Sscr(D) \setminus \bigcup_{\sigma \in D \setminus D'} \C^{\sigma^c} \times 0^{\sigma} \label{var_deg}
\end{align}
\end{Lemme}

\begin{proof}
The equality \eqref{var_deg_0} comes from the definition of $\Sscr$. We deduce the equality \eqref{var_deg} from the \eqref{var_deg_0} by using the stabily by $D'$ by taking the faces and hence, if $\sigma \notin D'$ then every cone containing $\sigma$ is not in $D'$. Then 
\[
\bigcup_{\sigma \in D \setminus D'} (\C^*)^{\sigma^c} \times 0^{\sigma}=\bigcup_{\sigma \in D \setminus D'} \bigcup_{\tau \succeq \sigma} (\C^*)^{\tau^c} \times 0^{\tau}=\bigcup_{\sigma \in D \setminus D'} \C^{\sigma^c} \times 0^{\sigma}
\]
\end{proof}

\begin{Rem}
The quantum toric stack obtained as quotient of $\Sscr(\Delta')$ can not have the same properties as the quantum toric stack described in $\Omega(D)$. In particular, if a generator of a cone tends to 0 then we lost the compactness (since the fans are not complete any more) or if a 1-cone of a cone tends to an other one cone then we get an non-separated stack (cf. \cite[Corollary 6.5 \& Example 6.20]{Toricstacks}). We will example of this in \cref{5-Ex}.
\end{Rem}

Let $h \in \overline{\Omega_I} \setminus \Omega_I$ and $D'(h)$ the maximal degenerated combinatorial type such that $h$ is $D'$-admissible (i.e. the degenerated combinatorial type from $D$ where we have removed every non strongly convex cones).

\begin{Lemme} \label{calcul_D'}
We have the following equality :
\[
D'(h)=D \setminus \{\sigma \in D \mid \dim\mathrm{Vect}(h(e_i), i \in \sigma)<|\sigma| \}
\]
\end{Lemme}
\begin{proof}
The combinatorial type $D'$ is obtained from $D$ by removing the elements $I$ of $D$ where $\Cone(h(e_i), i \in I)$ is not strongly convex. Since $D$ is simplicial, this condition is equivalent to the fact that the dimension of the space generated by the $h(e_i), i \in I$ is less than the cardinal of $I$. 
\end{proof}

The fiber of the universal family over the point $h$ has to be a quotient of $\Sscr(D'(h))$. Therefore we define  : \\
\begin{equation}
E_I\coloneqq \{(h,z) \in \overline{\Omega} \times \Sscr \mid z \in \Sscr(D'(h) \}
\end{equation}
and
\begin{equation}
    E_{IJ} \coloneqq E_I \cap (U_J \times \Sscr).
\end{equation}
By using the equalities \eqref{var_deg} and the one of \cref{calcul_D'}, we get 
\begin{Lemme}
$E_I$ can be be written as follows
\[
E_I=\{(h,z) \mid \forall \sigma \in D, \dim\mathrm{Vect}(h(e_i), i \in \sigma)<|\sigma| \Rightarrow (z_i, i \in \sigma) \neq 0  \} 
\]
\end{Lemme}

We define an action of $\C^{n-d}$ on $E_I \subset U_I \times \Sscr$ by :
\begin{equation} \label{action_EI}
t \cdot ([k],z)=([k],E(s_{\textrm{Mono},I}([k])(t))z)
\end{equation}
The projection $\pi : E_I \to \overline{\Omega_I}$ is invariant and hence descend to the quotient :
\[
\pcal_I : \overline{\Xscr_I}\coloneqq [E_I/\C^{n-d}] \to \overline{\Omega_I}
\]

In order to glue the different families $\overline{\Xscr}_I \to \overline{\Omega}_I$, we will study the behavior of the change of index of the section : 
\begin{Lemme}
Let $[k] \in \Omega_I \cap U_J$. Then
\[
s_{\textrm{Mono},J}([k])=s_{\textrm{Mono},I}([k]) \circ K_{JI}([k])^{-1}
\]
where $K_{JI}([k])=((s_{\textrm{Mono},I}[k])_{j}(e_i))_{i \in [\!1,n-d]\!], j \in J} \in \GL_{n-d}(\R)$
\end{Lemme}

\begin{proof}
We prove this in the same manner as \cref{prep_equiv}
\end{proof}

Since the action of $\GL_{n-d}(\R)$ on $\mathrm{Epi}(\R^{n-d},\R^n)$ is free, the matrices $(K_{JI}([k]))$ satisfy cocycle condition:
\begin{Lemme}
For all subset $I,J,L \subset \{1,\ldots,n\}$ of cardinal $n-d$ and for all $[k] \in U_I \cap U_J \cap U_L$, 
\[
K_{IL}([k])=K_{IJ}([k])K_{JL}([k])
\]
\end{Lemme}

Thanks to \cref{descente_quotient}, we deduce the following statement: 

\begin{Prop} \label{cocycle_quotient}
The identity $E_{IJ} \subset E_I \to E_{IJ} \subset E_J$ descends as a stack isomorphism
\begin{equation}
\Phi_{I,J} \colon \Xscr_{IJ}\coloneqq \left[E_I \cap U_J \times \Sscr/\C^{n-d}\right] \to \Xscr_{JI} \coloneqq \left[E_J \cap U_I \times \Sscr/\C^{n-d}\right]
\end{equation}
Moreover, for all $I,J,L$ subsets of $[\![1,n]\!]$ of cardinal $n-d$, the compose $\Phi_{I,J}\Phi_{J,L}\Phi_{L,I}$ is $2$-isomorphic to the identity.
\end{Prop}

The cocycles $\Phi_{IJ}$ define transition maps between the $\overline{\Xscr_I}$. We can define the following gluing:
\[
\overline{\Xscr} \coloneqq \mathop{\colim}\limits_{I,J} \left(\begin{tikzcd}
	{\mathscr{X}_{IJ}} & {\overline{\mathscr{X}_{I}}} \\
	{\mathscr{X}_{JI}} & {\overline{\mathscr{X}_{J}}}
	\arrow[hook, from=1-1, to=1-2]
	\arrow["{\Phi_{IJ}}"', from=1-1, to=2-1]
	\arrow[hook, from=2-1, to=2-2]
\end{tikzcd} \right)
\]
Note $E$ the union of $E_I$ i.e. 
\begin{equation}
E=\{(h,z) \in \overline{\Omega} \times \Sscr \mid \forall \sigma \in D, \dim\Vect(h(e_i), i \in \sigma)<|\sigma| \Rightarrow (z_i, i \in \sigma) \neq 0  \}
\end{equation}
The morphisms $\pcal_I$ glue into a morphism 
\begin{equation}
\pcal : \overline{\Xscr} \to \overline{\Omega}
\end{equation}
The group $G$ acts on $E$ in the following way:
\begin{equation} \label{action_sur_E}
  \sigma \cdot ([k],z)=({\sigma} \cdot[ k],z_{\sigma})  
\end{equation}

\begin{Lemme}
The action of $\sigma \in G$ send $E_I$ on $E_{\sigma^{-1}(I)}$.
\end{Lemme} 

\begin{proof}
Let $k : \R^{n-d} \to \R^n$ be a linear monomorphism and $\sigma^{-1}(i) \in \sigma^{-1}(I)$. The $\sigma^{-1}(i)$th component of $\sigma \cdot k=P_\sigma^{-1} k$ satisfies
\[
\forall x \in \R^{n-d}, (\sigma \cdot [k])_{\sigma^{-1}(i)}(x)=x_{i}
\]
We deduce that if $[k] \in U_I$ then $\sigma \cdot [k] \in U_{\sigma^{-1}(I)}$.
\end{proof}

\begin{Prop} \label{passage_quotient_univ_fam}
The action \eqref{action_sur_E} descend in an action on $\overline{\Xscr}$ and the projection $\pcal : \overline{\Xscr} \to \overline{\Omega}$ is a $G$-equivariant map which descends as a stack morphism : 
\begin{equation} \label{compactification_famille_universelle}
    \overline{\Xscr}/G \to \overline{\Omega}/G \eqqcolon \overline{\Mscr}
\end{equation}
\end{Prop}

\begin{proof}
The second part of the statement is clear. We will now prove the descent to quotient. For this, we firstly prive that every $\sigma$ defines a stack morphism $\overline{\Xscr_I} \to \overline{\Xscr_{\sigma^{-1}(I)}}$ for all $I$ : \\
Let $\sigma \in G$, $t \in \C^{n-d}$ and $([k],z) \in E_I$. Then
\begin{align*}
\sigma \cdot (t \cdot ([k],z))&=\sigma \cdot ([k],) \\
    &=(\sigma \cdot [k], \sigma \cdot E(s_{\mathrm{Mono},I} [k](t)z))=(\sigma \cdot [k], E(s_{\mathrm{Mono},\sigma^{-1}(I)}(\sigma \cdot [k])(t)z_\sigma)) \\
    &=t \cdot (\sigma \cdot [k],z_\sigma)=t\cdot (\sigma \cdot ([k],z))
\end{align*}
The  morphism $\sigma \cdot $ is a $\C^{n-d}$-equivariant morphism and which descend to quotient. \\
Moreover, since this morphism sends $E_{IJ}$ on $E_{\sigma^{-1}(I)\sigma^{-1}(J)}$ for all $I$ and $J$ then the isomorphisms $\overline{\Xscr_I} \to \overline{\Xscr_{\sigma^{-1}(I)}}$ glue into an automorphism of $\overline{\Xscr}$. 
\end{proof}

\begin{Thm} \label{Compactification}
The compactification $\overline{\Mscr}$ of $\Mscr$ is the base space of a familly (here, a smooth map where every fiber has the same dimension)
\begin{equation}
    \overline{\Xscr}/G \to \overline{\Mscr} 
\end{equation}
of quantum toric stacks such that 
\begin{itemize}
    \item Over $\Mscr$, the fibers are the quantum toric stacks given by a fan of combinatorial type $D$ ; 
    \item Over $\overline{\Mscr} \setminus \Mscr$, the fibers are quantum toric stacks given by a fan of a degenerated combinatorial type of $D$.
\end{itemize}

\end{Thm}

\subsection{Examples}
\label{5-EX}
\subsubsection{Projective spaces}

In this section, we will consider the combinatorial type $S_d$ of the projective space $\P^d$.

\begin{Lemme} \label{calcul_Pn_Ombar}
The right hand side of the equality \eqref{eq_ombar} is equal to:
\[
C=\left\{[p_0:\ldots:p_d] \in \Gr(1,\R^{d+1})=\R\P^d \mid p_i \geq 0, \sum_{i=0}^d p_i=1  \right\}
\]
\end{Lemme}

\begin{proof}
The closed subset $F$ is given by the equations
\[
\forall i<j, p_ip_j \geq 0
\]
since the kernel of the morphism $h : \R^{d+1} \to \R^d$ associated to a $S_d$-admissible morphism $h$ is generated by $\begin{pmatrix}
-h(e_{d+1})\\ 1
\end{pmatrix} \in \R_{\geq 0}^{d+1}$. 
Therefore we have $C \subset F$. Let us prove the other sense of the inclusion  : \\
Let $[p] \in F$. Up to taking a representative with an opposite sign,  we can suppose that, thanks to the equation of $F$,  the $p_i$ are non-negative. Then $\sum_{i=0}^d p_i$ is positive and:
\[
[p]=\left[\frac{p_0}{\sum_{i=0}^d p_i}:\ldots:\frac{p_d}{\sum_{i=0}^d p_i}\right] \in C
\]

\end{proof}
By the \cref{calcul_Pn_Ombar} and \cref{cnx_eq_ombar}, we get

\begin{Prop}
The space $\overline{\Omega(S_d)}$ is equal to 
\[
\left\{[\lambda_0:\ldots:\lambda_d] \in \Gr(1,\R^{d+1})=\R\P^d \mid \lambda_i \geq 0, \sum_{i=0}^d \lambda_i=1  \right\}
\]
Hence, the closed space $\overline{\Omega(S_d)}$ is a $d$-simplex.
\end{Prop}

\begin{Prop}
We get the family $\Xscr_{d-1}/\mathfrak{S}_{d} \to \overline{\Mscr_{d-1}}$ by projecting the degenerated coordinates of the restriction $(\Xscr_d/\mathfrak{S}_{d+1})_{\mid \overline{\Mscr_d}\setminus \Mscr_d} \to \overline{\Mscr_d}\setminus \Mscr_d$.
\end{Prop}

\begin{proof}
A toric stack over ${\overline{\Mscr_d}\setminus \Mscr_d}$ is obtained by at least one degeneration i.e. is a quotient of an open subset of $\C^d \setminus \{0\} \times \C$ (up to the permutation of the coordinates, we can suppose that the degeneration happened on the last  coordinate). Hence we can project:
\[\begin{tikzcd}
	{(\Xscr_d/\Sfrak_{d+1})_{|\overline{\Mscr}_d \setminus \Mscr_d}} && {\Xscr_{d-1}/\Sfrak_d} \\
	& {\overline{\Mscr}_d \setminus \Mscr_d}
	\arrow[from=1-1, to=2-2]
	\arrow[from=1-3, to=2-2]
	\arrow[from=1-1, to=1-3]
\end{tikzcd}\]
\end{proof}

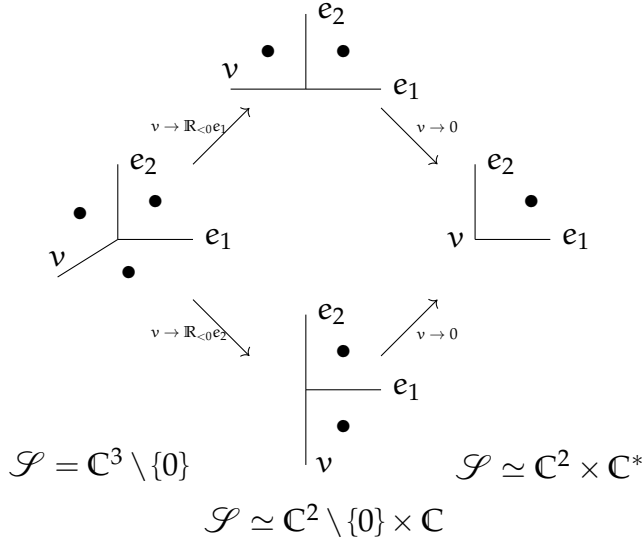
\begin{figure}[!h]
    \centering
\begin{tikzpicture}[scale=0.5]
\draw (0,0)  to (0,2) node[right]{$e_2$}  
node at (1,1) {$\bullet$}
;

   \draw (0,0) to (2,0) node[right]{$e_1$} 
   node at (-1,0.67) {$\bullet$}
;
   \draw (0,0)  to (-1.6,-1) node[right,above]{$v$} 
   node at (0.29,-.89){$\bullet$}
   ;
   
    \draw[->] (2, 2) to (3.5,3.5)
    node[scale=0.5] at (1.9,3) {$v \to \R_{<0}e_{1}$} ; 
    
   \draw[->] (2, -1.6) to (3.5,-3.1)
    node[scale=0.5] at (1.9,-2.5) {$v \to \R_{<0}e_{2}$} 
    node at (-0.5,-6) {$\Sscr=\C^3 \setminus \{0\}$} ; 
   
   \draw (5,4)  to (5,6) node[right]{$e_2$}  
   node at (6,5) {$\bullet$} ;
   \draw (5,4) to (7,4) node[right]{$e_1$} 
    node at (4,5) {$\bullet$} ;
   \draw (5,4)  to (3,4) node[right,above]{$v$} 
    node at (5.5,-7.5) {$\Sscr \simeq \C^2 \setminus \{0\} \times \C$} ; 
   
   \draw[->] (7,3.5) to (8.5, 2)
    node[scale=0.5] at (8.5,3) {$v \to 0$} ; 
   
   \draw (5,-4)  to (5,-2) node[right]{$e_2$}  
     node at (6,-3) {$\bullet$} ;
   \draw (5,-4) to (7,-4) node[right]{$e_1$} 
   node at (6,-5) {$\bullet$} ;
   \draw (5,-4)  to (5,-6) node[right]{$v$} 
    node at (11.5,-6) {$\Sscr \simeq \C^2 \times \C^*$} ; 
   
 \draw[->] (7,-3.1) to (8.5, -1.6)
    node[scale=0.5] at (8.5,-2.5) {$v \to 0$} ; 
    
   \draw (9.5,0)  to (9.5,2) node[right]{$e_2$}  
    node at (11,1) {$\bullet$} ;
   \draw (9.5,0) to (11.5,0) node[right]{$e_1$} 
   node[left] at (9.5,0) {$v$};
 \end{tikzpicture}  
\caption{Degeneration in dimension 2}
    \label{fig:my_label}
\end{figure}

\subsubsection{Example in dimension 2}
\label{5-Ex}
In this sub-subsection,we will suppose $d=2$. \\
We can remark that if we stay in the same chart, we have three possibilities of degeneration:
\begin{enumerate}
    \item vectors $h(e_i)$ becomes 0;
    \item two vectors $h(e_i)$ and $h(e_j)$, $\{i,j\} \in D$ become positively dependent;
    \item two vectors $h(e_i)$ and $h(e_j)$, $\{i,j\} \in D$ becomes negatively dependent.
\end{enumerate}
The three following propositions describe what can be said for each case: 

We will suppose, by simplicity, that the set of virtual generators is empty since they do not change the situation.

\begin{Prop}
Let $D$ be the combinatorial type of a complete fan of $\R^2$ and $D'$ be the degenerated combinatorial type obtained from $D$ when the vectors $h(e_i)$, $i \in I$ tends to 0. Then the obtained quantum toric stacks are no more complete and is of the form
\[
\Xscr \times \left[(\C^*)^I/\C^I \right]
\]
where $\C^I$ acts on $(\C^*)^I$ by the morphism $E$ and $\Xscr$ is described by the fan obtained by restricting the calibration on $\Z^{I^c} \oplus 0$ and by deleting the cones which contain the 1-cones generated by the $h(e_i), i \in I$.
\end{Prop}

\begin{proof}
Without loss of generality, we can suppose $I=\{n-|I|+1,\ldots,n\}$.
The calibration $h : \R^n \to \R^2$ is identically zero on $\R^I \oplus 0$. Consequently,  a basis of $\ker(h)$ is given by the concatenation of a basis of $h_{\R^{I^c}}$ with the vectors $e_i, i \in I$. The morphism of the Gale transform is the transpose of the map given by the matrix
\begin{equation}
 \begin{pmatrix}
v_1 \ldots v_{n-|I|-2} & 0 \\
 & Id_{|I|}
\end{pmatrix}   
\end{equation}
Then, the null vectors cannot be generators of $1$-cones so, by definition, the variety $\Sscr$ is of the form
\[
\Sscr' \times (\C^*)^I
\]
The action of $\C^{n-2}$ on $\Sscr$ can be decomposed in two actions: the first one on $\Sscr'$ by the action given by (the exponential of) the basis of $\ker(h_{|\R^{I^c}})$ and the second one on $(\C^*)^I$ given by (the exponential of) the identity. In other words,
\[
\Xscr_{h,\Delta}=[\Sscr/\C^{n-2}]=[\Sscr' \times (\C^*)^I/\C^{n-|I|-2} \times \C^{I}]=[\Sscr'/\C^{n-|I|-2}]\times [(\C^*)^I/\C^{I}]
\]
\end{proof}

\begin{Ex} \label{degen_0_4}
We will begin with
\[
D=\{1,2,3,4,(1,2),(2,3),(3,4),(4,1)\}
\]
When $h(e_3)$ tends to $0$, we get the combinatorial type
\[
D'=\{(1,2),(4,1)\}
\]
Then
\[
\Sscr(D')=\C_x \times \C_z^* \times \C_{y,t}^2 \setminus \{0\}.
\]

\[\begin{tikzpicture}
\draw (0,0)  to (0,2) node[right]{$e_2$}  ;
   \draw (0,0) to (2,0) node[right]{$e_1$} ;
   \draw (0,0)  to (-1.6,-1) node[right,above]{$v_4$} ;
   \draw[dashed] (0,0)  to (-.8,2) node[right,above]{$v_3 \to 0$} ;
\end{tikzpicture}
\]
\end{Ex}

\begin{Prop}

Let $D$ be the combinatorial type of a complete fan of $\R^2$ and $D'$ the degenerated combinatorial type obtained from $D$ when a couple $(h(e_i),h(e_j))$ tends to a couple of the form $(h(e_i), \lambda h(e_i))$, $\lambda >0$. 
Then we get a complete fan where the cone $\R_{\geq 0} h(e_i)$ appears twice and hence has two generators given by $h$ : $h(e_i)$ and $h(e_j)$.
\end{Prop}

\begin{Ex} \label{degen_+_4}
On part du type combinatoire 
\[
D=\{1,2,3,4,(1,2),(2,3),(3,4),(4,1)\}
\]
When $h(e_3)$ tends to the half-line $\R_{\geq 0} h(e_4)$, we get the combinatorial type
\[
D'=\{1,2,3,4,(1,2),(2,3),(4,1)\}
\]
\[\begin{tikzpicture}
\draw (0,0)  to (0,2) node[right]{$e_2$}  ;
   \draw (0,0) to (2,0) node[right]{$e_1$} ;
   \draw (0,0)  to (-1.6,-1) node[right,above]{$v_4$} node at (-0.8,-0.25) {$v_3$};
\end{tikzpicture}
\]

\end{Ex}

\begin{Rem}
We keep the two cones since it is prescribed by the combinatorics and the calibration 
\end{Rem}

\begin{Avert}
The obtained variety $\Sscr$ is different than the one obtained by keeping one generator. In the previous example, we have
\[
\Sscr(D')=\C^2 \times (\C^*)^2 \cup \C^* \times \C^2 \times \C^* \cup \C \times (\C^*)^2 \times \C
\]
and if we keep one generator (and by defining the other one as virtual), it is
\[
(\C^3 \setminus \{0\}) \times \C^*
\]
(the quotient is a quantum projective plane with one virtual generator). More precisely, we have the equality: 
\[
\Sscr(D') \setminus \{\C \times (\C^*)^2 \times 0 \}=[(\C^3 \setminus \{0\}) \times \C^* ]\setminus [0 \times \C^* \times 0 \times \C^*].
\]
\end{Avert}

\begin{Prop}
Let $D$ be the combinatorial type of a complete fan of $\R^2$ and $D'$ be the degenerated combinatorial type obtained from $D$ when a couple $(h(e_i),h(e_j))$ tends to a couple of the form $(h(e_i), \lambda h(e_i))$, $\lambda >0$.

If the group $\Gamma_{pr}\coloneqq\sum_{k \neq i,j} \left\langle w,h(e_k) \right\rangle\Z$ (where $w \perp h(e_i)$) is of rank 1 i.e. there exists $\alpha \in \Gamma_{pr}$, $\Gamma_{pr} =\alpha \Z$ then there exists a flat toric morphism $\Xscr_{\Delta',h} \to \Xscr_{\R_{\geq 0},\frac{1}{\alpha} \left\langle h,w\right\rangle}$ of generic fiber a quantum projective line and the fiber over $\left[0/\Z^{n-d}\right]$ is a finite union of quantum projective lines.
\end{Prop}

\begin{Rem}
\begin{itemize}
    \item This proposition is the irrationnal counterpart of the toric degeneration of \cite{torideg_Nishinou}.
    \item The condition on $\Gamma_{pr}$ is necessary in order to ensure that the morphism $\Xscr_{\Delta'} \to \Xscr_{\R_{\geq 0},\frac{1}{\alpha} \left\langle h,w\right\rangle}$ is well-defined (since the image of the generators of 1-cones in $ \Delta'$ have to be a integer multiple of a generator of $\R_{\geq 0}$.
\end{itemize}
\end{Rem}

\begin{proof}
Up to compose by an isomorphism, we can suppose $w=e_2$. The considered morphism  is the projection $\pi_2\colon \R^2 \to \R, (x,y) \mapsto y$. This is a fan morphism $(\Delta,h) \to (\R_{\geq 0},\pi_2 h : \Z^n \to \Gamma_{pr}$) if, and only if, for all $ i \in \Delta(1)$, $\pi_2 h(e_i) \in \N \pi_2 h(e_1)$. In this case, $ \Gamma$ is of rank 1. We can even suppose, without loss of generality, that is generated by 1.
In order to compute the fibre, we will use the following cartesian diagram :
\[
\begin{tikzcd}[row sep=scriptsize, column sep=scriptsize]
 	& \overline{\pi_2}^{-1}(x) \arrow[dl,hookrightarrow] \arrow[rr] \arrow[dd,"\overline{\pi_2}"' near start] & & \mathcal{p}^{-1}([x/\Z^{n-d}]) \arrow[dl,hookrightarrow] \arrow[dd,"\pcal"] \\
 	\mathscr{S} \arrow[rr, crossing over] \arrow[dd,"\overline{\pi_2}"'] & & {[\mathscr{S}/\C^{n-d}]}\\
 	& \Spec(\C) \arrow[dl,"{x}"'] \arrow[rr] & & \Bscr\Z^{n-d} \arrow[dl,hookrightarrow] \\
 	\C\arrow[rr] & & {[\C/\Z^{n-d}]} \arrow[from=uu, crossing over,"{\pcal}" near start]\\
 	\end{tikzcd}
 	\]
The morphism $\Sscr \to \C$ is flat (since it is non constant), we can deduce that the morphism on the quotient is also flat. \\
The fiber over a non-zero point is given by the cones $\R_{\geq 0} \times 0$ and $\R_{\leq 0} \times 0$. We have just to study the two cones containing them. Note $v=(v_1,v_2)$ the generator of the 1-cone which generates a cone with $\R_{\geq 0} \times 0$. By hypothesis, $v_2 \in \N_{>0}$ and thus the morphism $\overline{\pi_2}$ is defined by
\[
\overline{\pi_2}(z_1,z_2)=z_2^{v_2}
\]
So the fibre $\overline{\pi_2}^{-1}(E(t))$ is 
\begin{equation} \label{fibre}
   \C \times \mu_{v_2} E\left(\frac{t}{v_2}\right) 
\end{equation}
where $\mu_{v_2}$ is the group of $v_2$th root of the unity.
It remains to us to examine the action of $\Z^{n-d}$. It is given by the morphism
\[
\begin{pmatrix}
1 & v_1 \\ 0 & v_2
\end{pmatrix}^{-1}h=\begin{pmatrix}
1 & -v_1/v_2 \\ 0 & 1/v_2
\end{pmatrix}h
\]
This implies that the $v_2$ lines of \eqref{fibre} are identified by the action. We get the same result with the second cone. The transition maps glue these two lines in a quantum projective line. The case of the fiber over 0 is essentially the same. The only difference is that on the cone which is not in the border, the fiber over 0 is two secants lines.
\end{proof}
 \begin{Ex}  \label{degen_-_4}
 
We begin with the combinatorial type
\[
D=\{1,2,3,4,(1,2),(2,3),(3,4),(4,1)\}
\]
When $h(e_4)$ tends to the half-line $\R_{\leq 0} h(e_1)=\R_{\leq 0}e_1$, we get the combinatorial type
\[
D'=\{1,2,3,4,(1,2),(2,3),(3,4)\}
\]

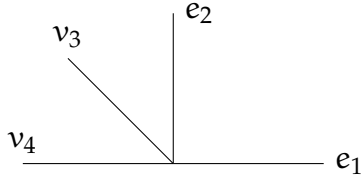
\begin{figure}[h!]
    \centering
\begin{tikzpicture}
\draw (0,0)  to (0,2) node[right]{$e_2$}  ;
   \draw (0,0) to (2,0) node[right]{$e_1$} ;
   \draw (0,0)  to (-2,0) node[right,above]{$v_4$} ;
   \draw (0,0)  to (-1.4,1.4) node[right,above]{$v_3$} ;
\end{tikzpicture}
\caption{Fan of combinatorial type $D'$}

\end{figure}

 The projection of the second coordinate is a toric morphism (as seen in the previous proposition). In order to find the fiber of this toric morphism, it suffice to décompose it in the following manner:
 
 \begin{figure}[h!]
     \centering

 \begin{tikzpicture}
    \draw (0,0)  to (0,1) node[right]{$e_2$}  ;
   \draw (0,0) to (1,0) node[right]{$e_1$} ;
   \draw (0,0)  to (-1,0) node[right,above]{$v_4$} ;
   \draw (0,0)  to (-0.7,0.7) node[right,above]{$v_3$} ;
   
    \draw[dotted,->] (1.5, 0.5) to (2.5,0.5);
    
    \draw (4,0) to (4,1) node[right]{$e_2$} ;
    \draw (4,0) to (5,0) node[right]{$e_1$} ;
    \draw (4,0) to (3,0) node[right,above]{$v_4$} ;
    
    \draw[->] (5.5, 0.5) to (6.5,0.5)
    node at (6,0.7) {$\pi$} ; 
    
     \draw (7,0) to (7,1) node[right]{$1$} ;
\end{tikzpicture}
\caption{Decomposition of the degeneration}
\end{figure}
 The first arrow is a blow-up of an equivariant map (when it exists) and the second one is the projection of the product. Outside zero, the fiber are quantum projective line and over 0, it is a quantum projective line where we have replaced an equivariant point (0 or $\infty$) by an other quantum projective line.
 
 \end{Ex}
 
 \subsubsection{Example with 4 generators}
 
 We will suppose $n=4$ and we will considere the following combinatorial type:
 \[C_4=\{\{1\},\{2\},\{3\},\{4\},\{1,2\},\{2,3\},\{3,4\},\{4,1\}\}
 \]
 This is the combinatorial type of fans of Hirzebruch surfaces. We have seen the degenerations of this combinatorial type in the examples \ref{degen_0_4},\ref{degen_+_4} and \ref{degen_-_4}. We will systematize this study in this paragraph. \\
 We have three cases to study (The other is obtained by action of $G$) : 
\begin{itemize}
    \item the usual chart \ $U_{(1,2)}$ ; 
    \item the chart $U_{(1,3)}$ ; 
    \item the chart $U_{(2,4)}$ ; 
\end{itemize}
The two last cases are treat in the same way (see figures \ref{fig:carte13} and \ref{fig:carte24}). Indeed, 
\[
\Omega_{(1,3)}=\Omega_{(2,4)}=\{((a,b),(c,d)) \in \R^4 \mid a,b>0,c,d<0 \}
\]
as illustrated by the figures \eqref{fig:carte13} and \eqref{fig:carte24}.
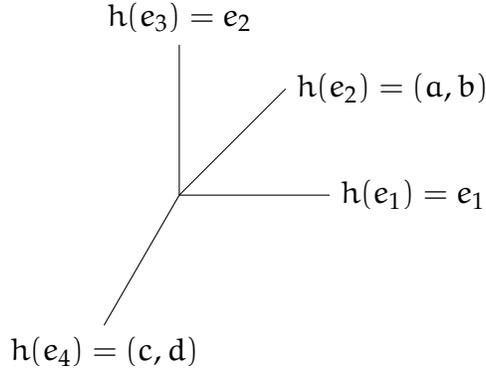
\begin{figure}[t!]
    \centering
    \begin{tikzpicture}
\draw (0,0)  to (0,2) node[above]{$h(e_3)=e_2$}  ;
   \draw (0,0) to (2,0) node[below,right]{$h(e_1)=e_1$} ;
   \draw (0,0)   to (1.414,1.414)node[below,right]{$h(e_2)=(a,b)$}  ;
   \draw (0,0)  to (-1,-1.731) node[below]{$h(e_4)=(c,d)$}  ;
\end{tikzpicture}
    \caption{Chart $U_{(1,3)}$}
    \label{fig:carte13}
\end{figure}
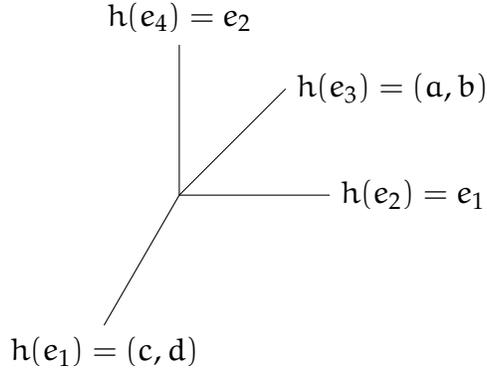
\begin{figure}[t!]
    \centering
    \begin{tikzpicture}
\draw (0,0)  to (0,2) node[above]{$h(e_4)=e_2$}  ;
   \draw (0,0) to (2,0) node[below,right]{$h(e_2)=e_1$} ;
   \draw (0,0)   to (1.414,1.414)node[below,right]{$h(e_3)=(a,b)$}  ;
   \draw (0,0)  to (-1,-1.731) node[below]{$h(e_1)=(c,d)$}  ;
\end{tikzpicture}
    \caption{Chart $U_{(2,4)}$}
    \label{fig:carte24}
\end{figure}
We deduce that
\begin{equation} \label{eq_13-24}
\overline{\Omega}_{(1,3)}=\overline{\Omega}_{(2,4)}=\R_{\geq 0}^2 \times \R_{\leq 0}^2.
\end{equation}
We will now focus on $\overline{\Omega}_{1,2}$. 
We know that $\overline{\Omega}_{1,2}$ is contained in the closed subset $F$ of $\R^{4}_{a,b,c,d}$ defined by $\begin{cases}
f_1(a,b,c,d)\coloneqq-a \geq 0 \\
f_2(a,b,c,d)\coloneqq ad-bc \geq 0 \\
f_3(a,b,c,d)\coloneqq-d \geq 0
\end{cases}$\\
 We can decompose $F$ in the following way:
 \begin{align*}
     F&= \coprod_{(\varepsilon_1,\varepsilon_2,\varepsilon_3) \in \{0,1\}^3 } \{(a,b,c,d) \mid \forall j, \mathrm{sign}(f_j(a,b,c,d))=\varepsilon_j \}   \\
     &= \Omega(D) \coprod \coprod_{(\varepsilon_1,\varepsilon_2,\varepsilon_3) \in \{0,1\}^3 \setminus \{(1,1,1) \}} \{(a,b,c,d) \mid \forall j, \mathrm{sign}(f_j(a,b,c,d))=\varepsilon_j \} 
 \end{align*}

It defines a stratification of $F$. 

\begin{Not}
We will note the set $\{(a,b,c,d) \mid \forall j, \mathrm{sign}(f_j(a,b,c,d))=\varepsilon_j \} $ by $F_{4\varepsilon_1+2\varepsilon_2+\varepsilon_3}$ (in other words, the binary representation of $k$ ($k \in\{0,\ldots,6\}$) gives the sign of the $f_i$ for $F_k$).
\end{Not}

\begin{Rem}
If $k=4\varepsilon_1+2\varepsilon_2+\varepsilon_3$ then the quantum toric stack associated $h \in F_i$ have $\varepsilon_1+\varepsilon_2+\varepsilon_3+1$ cones (the $+1$ come from the cone $\Cone(e_1,e_2)$ which is always here). More precisely, we remove the cones which corresponds to $i$ where $\varepsilon_i=0$
\end{Rem}
 \newpage

\begin{tabular}{|c|c|c|}
\hline \raisebox{-1.5cm}{$ F_0$} & \raisebox{-1.5cm}{$a=0$,$bc=0$,$d=0$} & 
\quad $b=0$, $c \leq 0$
    \hspace{-2cm}
    \raisebox{0.5cm}{\begin{tikzpicture}[scale=0.7]
\draw (0,0)  to (0,2) node[right]{$e_2$}  ;
   \draw (0,0) to (2,0) node[right]{$e_1$} ;
   \draw (0,0) node[below]{$v_3$}  to (-2,0) node[right,above]{$v_4$} ;
\end{tikzpicture}}

\hspace{1cm} $b=0$, $c \geq 0$
\!\!\!\!\!\!\!\!\!\!\!\!\!\!\!\!\!\!\raisebox{0.5cm}{\begin{tikzpicture}[scale=0.7]
\draw (0,0)  to (0,2) node[right]{$e_2$}  ;
   \draw (0,0) node[below,left]{$v_3$} to (2,0) node[right]{$e_1$} ;
   \draw (0,0)  to (1,0) node[right,above]{$v_4$} ;
\end{tikzpicture}}

\\
&&\quad $c=0$, $b \leq 0$
    \hspace{-2cm}
    \raisebox{0.5cm}{\begin{tikzpicture}[scale=0.7]
\draw (0,0) node[below,left]{$v_4$}  to (0,2) node[right]{$e_2$}  ;
   \draw (0,0) to (2,0) node[right]{$e_1$} ;
   \draw (0,0)  to (0,-2) node[right]{$v_3$} ;
\end{tikzpicture}}

\hspace{1cm} $c=0$, $b \geq 0$
\!\!\!\!\!\!\!\!\!\!\!\!\!\!\!\!\!\!\raisebox{0.5cm}{\begin{tikzpicture}[scale=0.7]
\draw (0,0)  to (0,2) node[right]{$e_2$}  ;
   \draw (0,0) node[below,left]{$v_4$} to (2,0) node[right]{$e_1$} ;
   \draw (0,0)  to (1,0) node[right,above]{$v_3$} ;
\end{tikzpicture}}
 \\ 
    \hline \raisebox{-1cm}{$ F_1$} & \raisebox{-1cm}{$a=0$,$bc=0$,$d<0$}  & \quad $b=0$
    \hspace{-2cm}
    \raisebox{0.5cm}{\begin{tikzpicture}[scale=0.7]
\draw (0,0)  to (0,2) node[right]{$e_2$}  ;
   \draw (0,0) to (2,0) node[right]{$e_1$} ;
   \draw (0,0) node[below, left]{$v_3$}  to (1,-1) node[right,above]{$v_4$} ;
\end{tikzpicture}}

\hspace{1cm} $b > 0, c=0$
\!\!\!\!\!\!\!\!\!\!\!\!\!\!\!\!\!\!\raisebox{0.5cm}{\begin{tikzpicture}[scale=0.7]
\draw (0,0)  to (0,2) node[right]{$e_2$}  ;
   \draw (0,0)  to (2,0) node[right]{$e_1$} ;
   \draw (0,0)  to (0,1) node[right]{$v_3$} ;
   \draw (0,0)  to (0,-1) node[right]{$v_4$} ;
\end{tikzpicture}}
 \\ 
 && \quad $b<0,c=0$
    \hspace{-2cm}
    \raisebox{0.5cm}{\begin{tikzpicture}[scale=0.7]
\draw (0,0)  to (0,2) node[right]{$e_2$}  ;
   \draw (0,0) to (2,0) node[right]{$e_1$} ;
   \draw (0,0) to (0,-1) node[right]{$v_3$} ;
   \draw (0,0) to (0,-2) node[left]{$v_4$} ;
\end{tikzpicture}}
 \\

 \hline
   \raisebox{1.7cm}{$F_2$}  & \raisebox{1.7cm}{$a=0$,$bc<0$,$d=0$}& \quad $b<0,c>0$
    \hspace{-2cm}
    \raisebox{0.5cm}{\begin{tikzpicture}[scale=0.7]
\draw (0,0)  to (0,2) node[right]{$e_2$}  ;
   \draw (0,0) to (2,0) node[right]{$e_1$} ;
   \draw (0,0) to (1,0) node[right,below]{$v_4$} ;
   \draw (0,0) to (0,-2) node[right]{$v_3$} ;
\end{tikzpicture}}

\hspace{1cm} $b > 0, c<0$
\!\!\!\!\!\!\!\!\!\!\!\!\!\!\!\!\!\!\raisebox{0.5cm}{\begin{tikzpicture}[scale=0.7]
\draw (0,0)  to (0,2) node[right]{$e_2$}  ;
   \draw (0,0)  to (2,0) node[right]{$e_1$} ;
   \draw (0,0)  to (0,1) node[right]{$v_3$} ;
   \draw (0,0)  to (-2,0) node[right,above]{$v_4$} ;
\end{tikzpicture}}
 \\  \hline
   
\end{tabular}

\begin{tabular}{|c|c|c|}
\hline
 \raisebox{1.7cm}{$F_3$} & \raisebox{1.7cm}{$a=0$,$bc<0$,$d<0$}&  \quad $b<0,c>0$
    \hspace{-2cm}
    \raisebox{0.5cm}{\begin{tikzpicture}[scale=0.7]
\draw (0,0)  to (0,2) node[right]{$e_2$}  ;
   \draw (0,0) to (2,0) node[right]{$e_1$} ;
   \draw (0,0) to (1,-1) node[right,below]{$v_4$} ;
   \draw (0,0) to (0,-2) node[right]{$v_3$} ;
\end{tikzpicture}}

\hspace{1cm} $b > 0, c<0$
\!\!\!\!\!\!\!\!\!\!\!\!\!\!\!\!\!\!\raisebox{0.5cm}{\begin{tikzpicture}[scale=0.7]
\draw (0,0)  to (0,2) node[right]{$e_2$}  ;
   \draw (0,0)  to (2,0) node[right]{$e_1$} ;
   \draw (0,0)  to (0,1) node[right]{$v_3$} ;
   \draw (0,0)  to (-1,-1) node[right,above]{$v_4$} ;
\end{tikzpicture}}
 \\   \hline
    \raisebox{-1cm}{$F_4$} & \raisebox{-1cm}{$a<0$,$bc=0$,$d=0$}&  \quad $b=0,c \geq 0$
    \hspace{-2cm}
    \raisebox{0.5cm}{\begin{tikzpicture}[scale=0.7]
\draw (0,0)  to (0,2) node[right]{$e_2$}  ;
   \draw (0,0) to (2,0) node[right]{$e_1$} ;
   \draw (0,0) to (1,0) node[right,below]{$v_4$} ;
   \draw (0,0) to (-2,0) node[right,below]{$v_3$} ;
\end{tikzpicture}}

\hspace{1cm} $b = 0, c \leq 0$
\!\!\!\!\!\!\!\!\!\!\!\!\!\!\!\!\!\!\raisebox{0.5cm}{\begin{tikzpicture}[scale=0.7]
\draw (0,0)  to (0,2) node[right]{$e_2$}  ;
   \draw (0,0)  to (2,0) node[right]{$e_1$} ;
   \draw (0,0)  to (-2,0) node[right,above]{$v_3$} ;
   \draw (0,0)  to (-1,0) node[right,above]{$v_4$} ;
\end{tikzpicture}}
 \\   
 && \quad $c=0, b \leq 0$
    \hspace{-2cm}
    \raisebox{0.5cm}{\begin{tikzpicture}[scale=0.7]
\draw (0,0)  to (0,2) node[right]{$e_2$}  ;
   \draw (0,0) to (2,0) node[right]{$e_1$} ;
   \draw (0,0) to (0,0) node[left]{$v_4$} ;
   \draw (0,0) to (-1,-1) node[right]{$v_3$} ;
\end{tikzpicture}} 
\hspace{1cm} $c = 0, b \geq 0$
\!\!\!\!\!\!\!\!\!\!\!\!\!\!\!\!\!\!\raisebox{0.5cm}{\begin{tikzpicture}[scale=0.7]
\draw (0,0)  to (0,2) node[right]{$e_2$}  ;
   \draw (0,0)  to (2,0) node[right]{$e_1$} ;
   \draw (0,0)  to (-1,1) node[right,above]{$v_3$} ;
   \draw (0,0)  to (0,0) node[left]{$v_4$} ;
\end{tikzpicture}}

\\
 \hline
    \raisebox{1.5cm}{$F_5$} & \raisebox{1.5cm}{$a<0$,$ad-bc=0$,$d<0$}& \quad $b>0,c>0$
    \hspace{-2cm}
    \raisebox{0.5cm}{\begin{tikzpicture}[scale=0.7]
\draw (0,0)  to (0,2) node[right]{$e_2$}  ;
   \draw (0,0) to (2,0) node[right]{$e_1$} ;
   \draw (0,0) to (2,-2) node[right,below]{$v_4$} ;
   \draw (0,0) to (-2,2) node[right,below]{$v_3$} ;
\end{tikzpicture}}

\hspace{1cm} $b <0, c < 0$
\!\!\!\!\!\!\!\!\!\!\!\!\!\!\!\!\!\!\raisebox{0.5cm}{\begin{tikzpicture}[scale=0.7]
\draw (0,0)  to (0,2) node[right]{$e_2$}  ;
   \draw (0,0)  to (2,0) node[right]{$e_1$} ;
   \draw (0,0)  to (-1,-1) node[right,above]{$v_3$} ;
   \draw (0,0)  to (-2,-2) node[right,above]{$v_4$} ;
\end{tikzpicture}}
  
 \\ 

 \hline
    \raisebox{1.2cm}{$F_6$} & \raisebox{1.2cm}{$a<0$,$bc<0$,$d=0$}& \quad $b>0,c<0$
    \hspace{-2cm}
    \raisebox{0.5cm}{\begin{tikzpicture}[scale=0.7]
\draw (0,0)  to (0,2) node[right]{$e_2$}  ;
   \draw (0,0) to (2,0) node[right]{$e_1$} ;
   \draw (0,0) to (-2,0) node[right,below]{$v_4$} ;
   \draw (0,0) to (-2,2) node[right,below]{$v_3$} ;
\end{tikzpicture}}

\hspace{1cm} $b <0, c > 0$
\!\!\!\!\!\!\!\!\!\!\!\!\!\!\!\!\!\!\raisebox{0.5cm}{\begin{tikzpicture}[scale=0.7]
\draw (0,0)  to (0,2) node[right]{$e_2$}  ;
   \draw (0,0)  to (2,0) node[right]{$e_1$} ;
   \draw (0,0)  to (-1,-1) node[right,above]{$v_3$} ;
   \draw (0,0)  to (1,0) node[right,above]{$v_4$} ;
\end{tikzpicture}}
 \\  \hline
\end{tabular}
\\

We deduce:

\begin{Prop} \label{adh12}
The closure $\overline{\Omega}_{1,2}$ is equal to the closed subset $F$ of the beginning of the sub-subsection.
\end{Prop}

The equalities \eqref{eq_13-24} and \cref{adh12} allow us to get the following equality: 
\begin{Thm}
The compactification $\overline{\Omega}$ is equal to
\[
\{ [a:b:c:d:e:f] \in \Gr(2,4) \subset \R\P^5 \mid ef\geq 0, af \leq 0, df\leq 0, ad \geq 0,ae \leq 0,de \leq 0\}
\]
\end{Thm}

\phantomsection{~}

\appendix
\section{Group action on quotient stacks}
\label{5-prelim}
This appendix will be dedicated to the study of group action on quotient stack which is not equivariant but verified a weaker assumption of compatibility (as in \cref{etape_calcul}). 

We will begin with a general statement on descent to quotient of smooth map.

\begin{Thm} \label{descente_quotient}
Let $G$ be a Lie group acting on a manifold $X$. Let $f : X \to X$ be a smooth map such that there exists a smooth map $K : X \to \Aut(G)$ such that 
\begin{equation} \label{presque_equiv}
    \forall x \in X, \forall g \in G, f(g \cdot x)=K(x)(g)f(x).
\end{equation}
Then $f$ descend in a stack morphism $[X/G] \to [X/G]$.
\end{Thm}

If $K$ is constant equal to the identity morphism i.e. $f$ is equivariant, it is immediate : to each object over $S$ of $[X/G]$
\[\begin{tikzcd}[ampersand replacement=\&]
	P \& X  \\
	S
	\arrow[from=1-1, to=2-1]
	\arrow["{m_P}", from=1-1, to=1-2]
\end{tikzcd}\]
we associate 
\[\begin{tikzcd}[ampersand replacement=\&]
	P \& X \& Y \\
	S
	\arrow[from=1-1, to=2-1]
	\arrow["{m_P}", from=1-1, to=1-2]
	\arrow["f", from=1-2, to=1-3]
\end{tikzcd}\]
(since $f$ is equivariant, the map $fm_P$ is equivariant). We do the same thing for morphism and we get a stack morphism $[X/G] \to [X/G]$ which verifies
\[\begin{tikzcd}[ampersand replacement=\&]
	X \& X \\
	{[X/G]} \& {[X/G]}
	\arrow["f", from=1-1, to=1-2]
	\arrow[from=1-1, to=2-1]
	\arrow[from=1-2, to=2-2]
	\arrow["{\overline{f}}"', from=2-1, to=2-2]
\end{tikzcd}\]

\begin{proof}

Firsty, we recall the universal property of the quotient of stacks by a group action (cf. \cite[Proposition 2.6]{Romagny}) i.e. the equivalence of categories 
\[\Hom_{G-\mathbf{Stack}}(X,[X/G]^{triv}) \simeq \Hom_{\mathbf{Stack}}([X/G],[X/G])
 \]
where $G-\mathbf{Stack}$ is the 
2-category of stacks with an action of $G$ and $[X/G]^{triv}$ is $[X/G]$ endowed with a trivial $G$-action.
Thank to this, it is enough to prove that the map $\pi f : X \to [X/G]$ (where $\pi$ is the projection $X \to [X/G]$ ) is $G$-invariant i.e.
we have the following 2-commutative diagram :

\begin{equation} \label{diagram_inv}
\begin{tikzcd}
	{X \times G} & X \\
	{[X/G] \times G} & {[X/G]}
	\arrow["{act}",from=1-1, to=1-2]
	\arrow[""{name=0, anchor=center, inner sep=0}, "{\pi f}", from=1-2, to=2-2]
	\arrow["{(\pi f,id)}"', from=1-1, to=2-1]
	\arrow[""{name=1, anchor=center, inner sep=0}, "{pr_1}"', from=2-1, to=2-2]
	\arrow[shorten <=4pt, shorten >=4pt, Rightarrow, from=0, to=1]
\end{tikzcd}
\end{equation}
By Yoneda, it is enough to prove that the fiber bundle $(X \times G) {}_{\pi f \mathrm{act}}\times_\pi X \to X \times G$ and  $(X \times G) {}_{(\pi f,id) act}\times_{\pi} X \to X \times G$ are isomorphic in a compatible way with the projection on $X$. \\
The first step of the computation of these fiber products is the computation of $X {}_{\pi f }\times_\pi X$ : 
\begin{equation} \label{produitfibre_pf}
\begin{tikzcd}
	{X \times G} & {} & {X \times G} & X \\
	X && X & {[X/G]}
	\arrow["\pi", from=1-4, to=2-4]
	\arrow["\pi", from=2-3, to=2-4]
	\arrow["act", from=1-3, to=1-4]
	\arrow["{pr_1}"', from=1-3, to=2-3]
	\arrow["f"', from=2-1, to=2-3]
	\arrow["{f \times id}", from=1-1, to=1-3]
	\arrow["{pr_1}"', from=1-1, to=2-1]
	\arrow["\lrcorner"{anchor=center, pos=0.125}, draw=none, from=1-3, to=2-4]
	\arrow["\lrcorner"{anchor=center, pos=0.125}, draw=none, from=1-1, to=2-3]
\end{tikzcd}
\end{equation}
\\
The fiber product on the right is the classical result used in order to give a presentation by a Lie groupoid of the quotient stack $[X/G]$ and the second one is the pullback of a trivial fiber bundle by the map $act$. In what follows, we will note $\alpha$ the map defined by
\[
\forall g \in G, \forall x \in X, \alpha(g,x)=g \cdot f(x)
\]
\\
We can compute the desired fiber product :
\\
\[\begin{tikzcd}
	{X \times G \times G} & {} & {X \times G} & X \\
	{X \times G} && X & {[X/G]}
	\arrow["\pi", from=1-4, to=2-4]
	\arrow["{\pi f}"', from=2-3, to=2-4]
	\arrow["\alpha", from=1-3, to=1-4]
	\arrow["{pr_1}"', from=1-3, to=2-3]
	\arrow[""{name=0, anchor=center, inner sep=0}, "act"', from=2-1, to=2-3]
	\arrow["\lrcorner"{anchor=center, pos=0.125}, draw=none, from=1-3, to=2-4]
	\arrow["{act \times id}", from=1-1, to=1-3]
	\arrow["{pr_{1,2}}"', from=1-1, to=2-1]
	\arrow["\lrcorner"{anchor=center, pos=0.125}, draw=none, from=1-1, to=0]
\end{tikzcd}\]
(The right fiber product is given by \eqref{produitfibre_pf} and the left one is the pullback of the trivial fiber bundle by $f$) and
\[\begin{tikzcd}
	{X \times G\times G} && {X \times G} & X \\
	{X \times G} && {[X/G] \times G} & {[X/G]}
	\arrow["\pi", from=1-4, to=2-4]
	\arrow[""{name=0, anchor=center, inner sep=0}, "{\pi f\times id}"', from=2-1, to=2-3]
	\arrow["{\pi \times id}", from=1-3, to=2-3]
	\arrow["{pr_1}"',from=2-3, to=2-4]
	\arrow["{pr_1}",from=1-3, to=1-4]
	\arrow["{\alpha \times id}", from=1-1, to=1-3]
	\arrow["{pr_{1,3}}"', from=1-1, to=2-1]
	\arrow["\lrcorner"{anchor=center, pos=0.125}, draw=none, from=1-3, to=2-4]
	\arrow["\lrcorner"{anchor=center, pos=0.125}, draw=none, from=1-1, to=0]
\end{tikzcd}\]
The right fiber bundle is the pullback of the trivial fiber bundle by $\pi$ and the left one is given by \eqref{produitfibre_pf}.
\\
Note $m_1$ the map $X \times G \times G \to X$ defined by
\[
(x,g,h) \mapsto \alpha \circ (act \times id)(x,g,h)=h \cdot f(g \cdot x) 
\]
and $m_2$ the map $X \times G \times G \to X$ defined by
\[
(x,g,h) \mapsto pr_1 \circ (\alpha \times id)(x,g,h)=g\cdot f(x).
\]
Note $\Phi$ the map $X \times G \times G \to X \times G \times G$ defined by :
\[
\forall (x,g,h) \in X \times G \times G, \Phi(x,g,h)=(x,hK(x)(g), h).
\]
Then the diagram 
\[\begin{tikzcd}
	& X \\
	{X \times G \times G} && {X \times G \times G} \\
	& {X \times G}
	\arrow["{pr_{1,3}}", from=2-3, to=3-2]
	\arrow["{pr_{1,2}}"', from=2-1, to=3-2]
	\arrow["\Phi", from=2-1, to=2-3]
	\arrow["{m_1}", from=2-1, to=1-2]
	\arrow["{m_2}"', from=2-3, to=1-2]
	\arrow["\simeq"', from=2-1, to=2-3]
\end{tikzcd}\]
commute. Indeed, the bottom triangle is clearly commutative and in order to prove the commutativity of the second one, we use the following equalities :
\begin{align}
    \forall (x,g,h) \in X \times G \times G,
    m_2(\Phi(x,g,h))&=m_2(x,hK(x)(g),h) \\
    &=hK(x)g \cdot f(x) \label{eq4_equiv}\\ 
    &=hf(g \cdot x)=m_1(x) \label{eq5_equiv}
\end{align}
The passage from the equality \eqref{eq4_equiv} to the equality \eqref{eq5_equiv} is done by the equality \eqref{presque_equiv}. \\
We deduce that the diagram \eqref{diagram_inv} is 2-commutative. That concludes the proof.
\end{proof}

Now we will consider the case of this paper: 
Let $X$ be a $G$-manifold and $H$ a Lie group acting on $X$ such that for all $\gamma \in H$, there exists $K_\gamma : X \to \Aut(G)$ such that
\[
\forall x \in X,\forall \gamma \in H, \forall g \in G, \gamma \cdot (gx)=K_\gamma(x)(g)\gamma \cdot x
\]
Then for all pair $(\gamma,\gamma') \in H^2$, we have :
\begin{equation} \label{action_cocycle}
    \gamma'\circ \gamma(gx)=\gamma'(K_\gamma(x)(g)\cdot (\gamma \cdot x))=K_{\gamma'}(\gamma \cdot x)(K_\gamma(x)(g))\cdot (\gamma'\gamma \cdot x)
\end{equation}

We can remark that the right hand side depend of $\gamma$ and $\gamma'$ but not of their product. Thus, the descent to quotient of this group action is not necessarily a group action.

\begin{Prop} \label{descente_action}

Suppose
\begin{equation} \label{eq_2-comm}
    \forall x \in X, \forall g \in G,\forall \gamma,\gamma'\in H, K_{\gamma'}(x)(K_\gamma(x)(g))=K_{\gamma\gamma'}(x)(g)
\end{equation}
Then the descent to quotient of the action of $H$ on $X$ is a group action on $[X/G]$. 
\end{Prop}

\begin{proof}
The condition on the identity element is clear by passing to quotient the identity $e\cdot x=x$ (for all $x$) and the equalities $\eqref{action_cocycle}$, $\eqref{eq_2-comm}$ (and the proof of \cref{descente_quotient})  give isomorphisms
\[
\forall \gamma,\gamma' \in H, \forall x : S \to [X/G],\gamma'\cdot (\gamma \cdot x) \Rightarrow \gamma \gamma' \cdot x.
\]
\end{proof}

\bibliographystyle{alpha} 
\bibliography{Biblio}

\end{document}